\documentclass[11pt]{amsart}
\usepackage{amsmath,amssymb,amsthm,mathrsfs}
\usepackage[all]{xy}
\usepackage[dvips]{graphicx}

\newtheorem{Def}{Definition}[section]
\newtheorem{Thm}[Def]{Theorem}
\newtheorem{Rem}[Def]{Remark}
\newtheorem{Prop}[Def]{Proposition}
\newtheorem{Conj}[Def]{Conjecture}
\newtheorem{Ex}[Def]{Example}

\newcommand{\C}{\mathbb{C}}
\newcommand{\D}{\mathbb{D}}
\newcommand{\R}{\mathbb{R}}

\newcommand{\HH}{\mathbb{H}}
\newcommand{\Z}{\mathbb{Z}}
\newcommand{\N}{\mathbb{N}}
\newcommand{\DD}{\mathbb{D}}

\newcommand{\PP}{\mathbb{P}}
\newcommand{\RR}{\mathbb{R}}
\newcommand{\Hol}{\mathrm{Hol}}

\newcommand{\End}{\mathrm{End}}
\newcommand{\Hom}{\mathrm{Hom}}
\newcommand{\Ker}{\mathrm{Ker}}
\newcommand{\U}{\mathrm{U}}

\begin{document}
\title[Degenerations, theta functions and geometric quantization]{Degenerations, theta functions and geometric quantization in mirror symmetry}
\author{Atsushi Kanazawa}

\maketitle

\begin{abstract}
We discuss various topics on degenerations and special Lagrangian torus fibrations of Calabi--Yau manifolds in the context of mirror symmetry. 
A particular emphasis is on Tyurin degenerations and the Doran--Harder--Thompson conjecture, which builds a bridge between mirror symmetry for Calabi--Yau manifolds and that for quasi-Fano manifolds.  
The proof of the conjecture is of interest in its own right and leads us to a few other related topics such as SYZ mirror symmetry, theta functions and geometric quantization.  
Inspired by the conjecture, we also propose a new construction of Landau--Ginzburg models by splitting Calabi--Yau fibrations.  
\end{abstract}

\tableofcontents

%%%%%%%%%%%%%%%%%%%%%%%%%%%%%%%%%%%%%%%%%%%%%%%%%%%%%%%%%%%%%%%%%%%%%%%%%%%%%%%%%%%%%%%%%%%%%
%%%%%%%%%%%%%%%%%%%%%%%%%%%%%%%%%%%%%%%%%%%%%%%%%%%%%%%%%%%%%%%%%%%%%%%%%%%%%%%%%%%%%%%%%%%%%

\section{Introduction}

We discuss various topics on degenerations and special Lagrangian torus fibrations of Calabi--Yau manifolds in the context of mirror symmetry. 
Mirror symmetry began as a phenomenon in superstring theory in late 1980s. 
Superstring theory posits that our spacetime is locally the product of the Minkowski space and a Calabi--Yau 3-fold. 
A surprising observation is that two superstring theories based on two distinct Calabi--Yau 3-folds sometimes gives rise to the same physical theory. 
A detailed study of this duality led us to the idea of mirror symmetry for Calabi--Yau manifolds, which can be loosely stated as follows.  
\begin{Conj}[Mirror Symmetry]
For a Calabi--Yau manifold $X$, there exists another Calabi--Yau manifold $Y$, called a mirror manifold, 
such that the complex geometry of $X$ is equivalent to the symplectic geometry of $Y$, and vice versa.
\end{Conj}
The above {\it equivalence} has been formulated and confirmed for many examples. 
A mirror pair of Calabi--Yau $n$-folds $X$ and $Y$ exhibits an exchange of Hodge numbers 
$$
h^{1,1}(X) = h^{n-1,1}(Y), \ \ \ h^{n-1,1}(X) = h^{1,1}(Y),
$$
and string theorists explicitly constructed many such pairs of Calabi--Yau 3-folds based on superconformal field theories \cite{GrePle, CLS}. 
This simple but elegant duality of the Hodge numbers immediately attracted much attention from mathematics. 
More surprisingly, in the celebrated work  \cite{CdOGP}, Candelas, de la Ossa, Green and Parkes computed the number of the rational curves of every fixed degree in a quintic Calabi--Yau 3-fold in $\PP^4$ 
by certain period integral calculations for the mirror Calabi--Yau 3-fold. 
Although the methods used in their work were based on physical intuition and thus not rigorous, they gave an amazing answer to a long-standing open problem in enumerative geometry. 
Their work greatly stunned algebraic and symplectic geometers, and there has followed more than two decades of very rewarding efforts to understand the mathematical mechanism underlying mirror symmetry.

There are various formulations of mirror symmetry and each one is interesting in its own right. 
The process of building a mathematical foundation of mirror symmetry has given impetus to new fields in mathematics, 
such as Gromov--Witten theory, Fukaya category and Bridgeland stability conditions. 
Mirror symmetry has also been a source of many new insights and progresses in algebraic geometry and symplectic geometry. 
In the development of mirror symmetry, it has also become more apparent that Calabi--Yau manifolds enjoy very rich properties. 

On the other hand, it has been noticed that mirror symmetry can be formulated for a much larger class of varieties such as Fano manifolds and varieties of general type.   
For example, the mirror of a Fano manifold with a choice of its anti-canonical divisor is given by a Landau--Ginzburg model, 
which is a pair consisting of a K\"ahler manifold $Y$ and a holomorphic function $W:Y\rightarrow \C$, called a superpotential \cite{EHX, Hor}. 
The study of Landau--Ginzburg models has reinvigorated many branches of mathematics such as singularity theory, matrix factorizations, and primitive forms. 

Today there are two principal approaches toward understanding the mechanism of mirror symmetry in mathematics. 
One is Kontsevich's homological mirror symmetry \cite{Kon} and the other is the Strominger--Yau--Zaslow (SYZ) mirror symmetry \cite{SYZ}.  
In this article we will mostly focus on the latter conjecture.  
\begin{Conj}[SYZ mirror symmetry \cite{SYZ}]
A Calabi--Yau $n$-fold $X$ admits a special Lagrangian $T^n$-fibration $\pi:X\rightarrow B$ 
and a mirror Calabi--Yau $n$-fold $Y$ is obtained as the total space of the dual $T^n$-fibration $\pi^\vee:Y\rightarrow B$. 
These fibrations are called SYZ fibrations. 
\end{Conj}
The heart of the SYZ conjecture is that mirror symmetry can be understood by dividing a Calabi--Yau $n$-fold into tori $T^n$ and then T-dualizing them to get the mirror Calabi--Yau $n$-fold.  
In this SYZ picture, mirror symmetry is thought of as a generalization of the Fourier transformation, relating various mathematical objects on distinct Calabi--Yau manifolds in highly non-trivial ways. 
The SYZ conjecture asserts that a Calabi--Yau manifold admits the structure of a completely integrable system (the Liouville--Arnold theorem), 
and moreover the level sets of the preserved quantities of the system are minimal submanifolds.  
Moreover, considering the fact that in a K3 surface the special Lagrangian tori and the elliptic curves are related by hyperK\"ahler rotations, 
we can regard the SYZ mirror symmetry as a vast structure theorem for Calabi--Yau manifolds in high dimensions, generalizing the study of elliptic K3 surfaces. 

Mirror symmetry is known to be intimately related to degenerations of Calabi--Yau manifolds.   
One important theme of this article is the SYZ fibrations can be approximated by certain degenerations of Calabi--Yau manifolds (Section \ref{Deg and SYZ}).  
Good examples are the Gross--Siebert program and the Doran--Harder--Thompson (DHT) conjecture. 
Since there is a good survey \cite{GroSie0} of the Gross--Siebert program by the founders, we will focus on the latter topic in this article. 
The basic setup is as follows. 
Given a Tyurin degeneration of a Calabi--Yau manifold $X$ to a union of two quasi-Fano manifolds $X_1 \cup_Z X_2$ 
intersecting along a common smooth anti-canonical divisor $Z \in |-K_{X_i}|$,  
it is natural to investigate a potential relationship between geometry of the Calabi--Yau manifold $X$ and that of the quasi-Fano manifolds $X_1$ and $X_2$.  
The DHT conjecture (Conjecture \ref{DHT}) builds a bridge between mirror symmetry for the Calabi--Yau manifold $X$ and that for the quasi-Fano manifolds $X_1$ and $X_2$ \cite{DHT}. 
It claims that we should be able to glue together the mirror Landau--Ginzburg models $W_i:Y_i \rightarrow \C$ of the pair $(X_i,Z)$ for $i=1,2$ 
to construct a mirror Calabi--Yau manifold $Y$ of $X$ equipped with a fibration $W:Y\rightarrow \PP^1$. 

The author recently proved the DHT conjecture for the elliptic curves by using ideas from SYZ mirror symmetry \cite{Kan1}. 
We will give a review of the proof and also extend it to the abelian surface case in this article. 
The key idea in the proof is twofold. 
\begin{enumerate}
\item The first is to obtain the correct complex structure by gluing the underlying affine base manifolds of $X_1$ and $X_2$ in SYZ mirror symmetry.  
This is based on the philosophy that a Tyurin degeneration of a Calabi--Yau manifold $X$ 
can be thought to be fibred over a Heegaard splitting of the base $B$ of a special Lagrangian torus fibration $\phi:X\rightarrow B$. 
\item The second is to construct theta functions out of the Landau--Ginzburg superpotentials. 
We observe that the product expressions of the theta functions are the manifestation of quantum corrections, which are encoded in the Landau--Ginzburg superpotentials, in SYZ mirror symmetry.
\end{enumerate}

An interesting feature of the proof is that theta functions are present in an unusual fashion. 
In fact, the appearance of theta functions in mirror symmetry is well-known (Section \ref{basis HMS}) 
and recent studies, especially Gross--Hacking--Keel--Siebert \cite{GHKS}, show that such canonical bases exist for a large class of Calabi--Yau manifolds. 
There is a nice survey on theta functions in mirror symmetry \cite{GroSie2}. 
However, we will take a quite different path, namely geometric quantization, to the theta functions.  
This circle of ideas was initially proposed in a series of works by Andrei Tyurin. 
The author thinks that geometric quantization is still lurking and only partially explored area in mirror symmetry.

\ \\
{\bf Structure of article} \\
In Section \ref{MS}, we recall some basics of mirror symmetry for both Calabi--Yau and Fano manifolds. To clarify the entire picture, we take a close look at K3 surfaces as working examples. 
In Section \ref{Deg and DHT} we discuss degenerations of Calabi--Yau manifolds and formulation of the DHT conjecture for Tyurin degenerations.  
In Section \ref{SYZ}, we gives a brief review of SYZ mirror symmetry, which will be a key tool in our proof of the DHT conjecture.  
In Section \ref{Deg and SYZ}, we provide a heuristic but important link between degenerations and SYZ fibrations of Calabi--Yau manifolds. 
In Section \ref{DHT proof}, we reconstruct the proof of the DHT conjecture, following \cite{Kan1} and extend it to the case of certain polarized abelian surfaces.  
In Section \ref{split GL}, we propose a new construction of Landau--Ginzburg models by using idea of the DHT conjecture. 
This works nicely for certain rational surfaces and possibly gives a powerful construction of Landau--Ginzburg models outside of the toric setting in higher dimensions
In Section \ref{GQ}, we discuss speculations on theta functions in geometric quantization and SYZ mirror symmetry. 
This section is mostly inspired by Tyurin's articles and is even more informal and speculative than other sections.   

Throughout the article, the author tries to keep precision combined with informality. 
Thereby some statements are not as formal as usually required and we prefer providing basic ideas and heuristics that are often hidden in research articles. 

\ \\
{\bf Acknowledgements} \\
The author would like to thank Shinobu Hosono, Hiroshi Iritani, Siu-Cheong Lau, Fumihiko Sanda, Shing-Tung Yau and Yuecheng Zhu for useful conversations on related topics. 
This survey is based on lectures delivered by the author at various places including the Harvard CMSA, KIAS, Kavli IPMU, Kyoto University, Kobe University, Gakushuin University and Hiroshima University. 
The author is grateful to these institutes for the warm hospitality and excellent research environment.  
This research is supported by the Kyoto Hakubi Project and JSPS Grant-in-Aid for Young Scientists(B) 17K17817.

%%%%%%%%%%%%%%%%%%%%%%%%%%%%%%%%%%%%%%%%%%%%%%%%%%%%%%%%%%%%%%%%%%%%%%%%%%%%%%%%%%%%%%%%%%%%%
%%%%%%%%%%%%%%%%%%%%%%%%%%%%%%%%%%%%%%%%%%%%%%%%%%%%%%%%%%%%%%%%%%%%%%%%%%%%%%%%%%%%%%%%%%%%%

\section{Mirror symmetry} \label{MS}

\subsection{Mirror symmetry for Calabi--Yau manifolds}
In this section we briefly review some basics of mirror symmetry for Calabi--Yau manifolds with particular emphasis on degenerations of Calabi--Yau manifolds.  
We refer the reader to \cite{CoxKat} for a detailed treatment of the subject. 

\begin{Def}
A Calabi--Yau manifold $X$ is a compact K\"ahler manifold such that the canonical bundle is trivial $K_X =0$. 
We sometimes assume $H^i(X,\mathcal{O}_X)=0$ for $0 < i <\dim X$. 
\end{Def}
Let $X$ be a Calabi--Yau $n$-fold.  
The triviality of the canonical bundle of $X$ implies that there exists a holomorphic volume form $\Omega$ up to multiplication by constants. 
It can be regarded as a complex orientation of $X$.  
We often think of a holomorphic volume form $\Omega$ and a K\"ahler form $\omega$ as part of the Calabi--Yau structure, especially when we discuss SYZ mirror symmetry.

There are two particularly interesting cohomology groups, namely $H^{1,1}(X)$ and $H^{n,1}(X)$. 
The former $H^{1,1}(X)$ represents the deformation of the complexified K\"ahler (symplectic) structure of $X$ 
since $H^{1,1}(X) \cap H^2(X,\R)$ contains the K\"ahler cone as an open cone, provided that $H^{2}(X,\mathcal{O})=0$. 
A complexified K\"ahler class $\omega$ is an element 
$$
\omega=B+\sqrt{-1}\kappa \in H^{1,1}(X)/2 \pi (H^{1,1}(X) \cap H^2(X,\Z)). 
$$
such that $\kappa$ is a K\"ahler class.  
The latter $H^{n,1}(X)$ is isomorphic to $H^1(X,T_X)$ by the Calabi--Yau condition $K_{X}= 0$ and represents the first order deformation of the complex structure of $X$. 
In fact, the Bogomolov--Tian--Todorov theorem asserts that the Kuranishi map
$$
K:H^1(X,T_X) \longrightarrow H^2(X,T_X)
$$
is the zero map. 
Thus the complex moduli space of a Calabi--Yau manifold $X$ is smooth and $H^1(X,T_X)$ represents the local complex moduli space around $X$. 

Mirror symmetry is a statement about Calabi--Yau manifolds in certain limits in the complex and K\"ahler moduli spaces. 
\begin{Conj}[Mirror Symmetry]
Given a Calabi--Yau manifold $X$ near a large complex structure limit, 
there exists another Calabi--Yau manifold $Y$, called a mirror manifold, such that 
complex geometry of $X$ is equivalent to symplectic geometry of $Y$ near the large volume limit, and vice versa. 
\end{Conj}

A large complex structure limit (LCSL) can be thought of as a point in the complex moduli space where $X$ maximally degenerates. 
A large volume limit (LVL) means a choice of a complexified K\"ahler class $\omega \in H^2(Y,\C)$ such that $\int_{C}\Im(\omega) \gg 0$ for every effective curve $C \subset Y$. 
We refer the reader to \cite{CoxKat} for more details of LCSLs and LVLs. 
In accordance with string theory, symplectic geometry and complex geometry are often called A-model and B-model.  
The simplest check of mirror symmetry is an exchange of the dimensions of the A- and B-model moduli spaces
$$
h^{1,1}(X) = h^{n-1,1}(Y), \ \ \ h^{n-1,1}(X) = h^{1,1}(Y).
$$ 
This is often called topological, or Hodge number mirror symmetry. 

\begin{Ex}[Batyrev \cite{Bat}] \label{Batyrev}
We recall some notation from toric geometry, which we will use in this article. 
Let $M \cong \Z^n$ be a lattice and $N=\Hom(M,\Z)$ the dual lattice. 
A lattice polytope $\Delta \subset M_\R=M\otimes \R$ is a convex hull of finitely many lattice points in $M$.
We define its polar dual by 
$$
\Delta^\vee=\{y \in \N_\R \ | \ \langle x,y\rangle \ge -1, \ \forall x \in \Delta\}.
$$ 
We say $\Delta$ is reflexive if $0 \in M_\R$ is in the interior of $\Delta$ and $\Delta^\vee$ is a lattice polytope. Note that $\Delta$ is reflexive if and only if $\Delta^\vee$ is. 

The toric variety $\PP_\Delta$ associated to a reflexive polytope $\Delta$ is a Fano variety. 
Batyrev proved that for a reflexive polytope $\Delta$ in dimensions $\le 4$ a general anti-canonical hypersurface $X'_\Delta \subset \PP_\Delta$ admits a crepant Calabi--Yau resolution $X_\Delta$. 
In the $3$-dimensional case, he confirmed the duality 
$$
h^{1,1}(X_\Delta) = h^{2,1}(X_{\Delta^\vee}), \ \ \  h^{2,1}(X_\Delta) = h^{1,1}(X_{\Delta^\vee}). 
$$
Moreover, Aspinwall--Greene--Morrison introduced a nice combinatorial correspondence between the complex and K\"ahker moduli spaces, called the monomial-divisor mirror map \cite{AGM}. 
Therefore mirror symmetry for this class of Calabi--Yau $3$-folds is at a combinatorial level elegantly simple. 
\end{Ex}

Surprisingly, the Batyrev construction produces 473,800,776 examples of mirror pairs of Calabi--Yau 3-folds, including 30,108 distinct pairs of Hodge numbers. 
In fact, it is an open problem whether or not the number of topological types of Calabi--Yau $3$-folds is bounded \cite{KanWil, Wil}. 
This Batyrev mirror construction was later generalized for the complete intersection Calabi--Yau manifolds in the toric Fano manifolds by Batyrev and Borisov \cite{BatBor}. 
The Batyrev--Borisov construction lay the mathematical foundation for much of the future research in mirror symmetry, and provide an excellent testing-ground for new conjectures and theories.  

Mirror symmetry should involve more than a mere exchange of Hodge numbers. 
A slightly refined version is the Hodge theoretic mirror symmetry, 
claiming an equivalence of the A-model Hogde structure associated to quantum cohomology of $Y$ 
and the B-model Hodge structure associated to the period integrals (the variation of Hodge structures) of $X$.  
This duality captures much finer information and involves the so-called mirror map 
which locally identifies the complex moduli space $\mathcal{M}_{\text{cpx}}(X)$ near a specified LCSL and the K\"ahler moduli spaces $\mathcal{M}_{\text{K\"ah}}(Y)$ near the LVL.  
For example, the famous calculation of Candelas and his collaborators can be understood in this framework \cite{CdOGP, Mor1}. 
The Hodge theoretic mirror symmetry, which in particular implies that the $g=0$ Gromov--Witten invariants of $X$ can be computed by certain period integrals of the mirror $Y$, 
is confirmed for a large class of Calabi--Yau manifolds by Givental \cite{Giv} and Lian--Liu--Yau \cite{LLY}. 

It is important to keep in mind that mirror correspondence depends upon the choice of a LCSL. 
Therefore mirror symmetry is inherently related to degenerations of Calabi--Yau manifolds. 
This also explains the failure of the conventional mirror symmetry for the rigid Calabi--Yau manifolds. 
If the complex moduli space $\mathcal{M}_{\text{cpx}}(X)$ of a Calabi--Yau manifold $X$ has several LCSLs, there should be several mirror manifolds $Y_1,Y_2,\dots$ accordingly. 
However, the existence of a LCSL is highly non-trivial, and in fact there exists a Calabi--Yau manifold whose (non-trivial) complex moduli space does not have such a point \cite{CvS}.  
On the other hand, there is a Calabi--Yau manifold whose complex moduli space has more than one LCSL. 
The first example was discovered by R\o dland \cite{Rod} and recently several more examples were constructed (see for example \cite{HosTak1,Kan1,Miu}).   

\begin{Ex}[Pfaffian--Grassmannian double mirror \cite{Rod}] \label{GP}
The Grassmannian $\mathrm{Gr}(2,7)$ of 2-dimensional subspaces in $\C^7$ has a canonical polarization via the Pl\"ucker embedding into $\PP(\wedge^2 \C^7)\cong \PP^{20}$. 
Let $L \subset \wedge^2 \C^7$ be a $13$-dimensional subspace of $\wedge^2 \C^7$ in a general position. 
Then 
$$
X_1=\mathrm{Gr}(2,7)\cap \PP(L) \subset \PP(\wedge^2 \C^7)
$$
is a Grassmannian Calabi--Yau 3-fold. 

On the other hand, the projective dual of $\mathrm{Gr}(2,7)$ in the dual space $\PP(\wedge^2 (\C^*)^7)$ is the Pfaffian variety $\mathrm{Pfaff}(7) \subset \PP(\wedge^2 (\C^*)^7)$. 
Another way to see $\mathrm{Pfaff}(7)$ is the rank $4$ locus of $\PP(\wedge^2 (\C^*)^7)$ when we identify $\wedge^2 (\C^*)^7$ with the space of skew-symmetric linear maps $\C^7 \rightarrow (\C^*)^7$. 
Let $L^\perp \subset \wedge^2 (\C^*)^7$ be the $13$-dimensional subspace dua to $L$. 
Then 
$$
X_2=\mathrm{Pfaff}(7)\cap \PP(L^\perp) \subset \PP(\wedge^2 (\C^*)^7)
$$
is a Pfaffian Calabi--Yau 3-fold. 

A relatively straightforward computation shows that the coincidence of Hodge numbers
$$
h^{1,1}(X_1)=h^{1,1}(X_2)=1, \ \ \ 
h^{2,1}(X_1)=h^{2,1}(X_2)=50.
$$ 
On the other hand, $H_1^3=42$ and $H_2^3=13$, where $H_1$ and $H_2$ are the ample generators of the Picard groups of $X_1$ and $X_2$ respectively. 
This implies that $X_1$ and $X_2$ are not birational to each other. 

Interestingly it is observed that the two Calabi--Yau 3-folds $X_1$ and $X_2$ share the same mirror Calabi--Yau 3-fold $Y$ whose complex moduli space has exactly two LCSLs: 
one corresponds to $X_1$ and the other to $X_2$. 
This double mirror phenomenon was mathematically confirmed in \cite{BCFKS,Tjo}. 
\end{Ex}

There are two principal approaches towards understanding the mechanism of mirror symmetry. 
One is Kontsevich's homological mirror symmetry \cite{Kon} and the other is the Strominger--Yau--Zaslow (SYZ) mirror symmetry \cite{SYZ}. 
Since we will provide a review on SYZ mirror symmetry in Section 4, we shall take a moment to look at homological mirror symmetry here. 

Calabi--Yau $n$-folds $X$ and $Y$ are called homological mirror symmetric if there is an equivalence of triangulated categories
$$
\mathrm{D^bCoh}(X)\cong \mathrm{D^bFuk}(Y). 
$$
Here $\mathrm{D^bCoh}(X)$ denotes the derived category of coherent sheaves on $X$ and $\mathrm{D^bFuk}(Y)$ denotes the derived Fukaya category of Lagrangian submanifolds in $Y$. 
We may think of the above equivalence as an categorification of an isomorphism 
$$
(K(X),\chi) \cong (H_n(Y,\Z),\cap)
$$
of lattices (modulo torsions), where $\chi$ is the Euler form and $\cap$ is the intersection pairing. 
Homological mirror symmetry is nothing but an equivalence of the B-branes and A-branes, 
but it is worth noting that Kontsevich conjectured such an equivalence before the discovery of the D-branes in superstring theory. 
It is very hard to find a homological mirror pair, and a fundamental work \cite{Sei3} of Seidel confirmed a version of homological mirror symmetry for a quartic hypersurface $Y \subset \PP^3$ and its mirror $X$.  
A version of homological mirror symmetry for a genus two curve was shown again by Seidel \cite{Sei2} and later generalized to the higher genus curves by Efimov \cite{Efi}.  

Assume that Calabi--Yau manifolds $X_1$ and $X_2$ share the same mirror manifold $Y$ (see for instance Example \ref{GP}). 
Then homological mirror symmetry implies that 
$$
\mathrm{D^bCoh}(X_1)\cong \mathrm{D^bFuk}(Y) \cong \mathrm{D^bCoh}(X_2),
$$
and thus $X_1$ and $X_2$ are Fourier--Mukai partners. 
Recall that non-isomorphic Calabi--Yau manifolds $X$ and $X'$ are said to be Fourier--Mukai partners if there is an equivalence 
$$
\mathrm{D^bCoh}(X)\cong \mathrm{D^bCoh}(X').
$$ 
We have a quite good understanding of the Fourier--Mukai partners of K3 surfaces via the lattice theory \cite{Ogu2,HLOY}, but very little is known about the Fourier--Mukai partners of Calabi--Yau manifolds in higher dimensions. 
One important result is that birational Calabi--Yau $3$-folds are Fourier--Mukai partners \cite{Bri0}.    
This is compatible with the well-known observation that mirror symmetry does not distinguish birational Calabi--Yau manifolds. 
%However, birational Calabi--Yau $3$-folds are not particularly interesting Fourier--Mukai partners from the view point of mirror symmetry.  

\begin{Ex}
It is shown by Borisov and C\u{a}ld\u{a}raru \cite{BorCal} that the Grassmannian Calabi--Yau 3-fold $X_1$ and the Pfaffian Calabi--Yau 3-fold $X_2$ are derived equivalent as mirror symmetry implied (Example \ref{GP}).  
Recall that they are not birational to each other. 
Later Hosono and Konishi calculated the higher genus Gromov--Witten invariants of $X_1$ and $X_2$ by solving the BCOV holomorphic equation \cite{BCOV, HosKon, KanZho}. 
They observed an interesting swapping of the higher genus Gromov--Witten invariants of the two Calabi--Yau 3-folds. 
\end{Ex}

Another interesting example was recently constructed by Hosono and Takagi \cite{HosTak1}. 
Their example is based on the projective duality between the secant varieties of symmetric forms and these of the dual forms. 
In this setting they naturally come into two Calabi--Yau 3-folds which are derived equivalent but not birational to each other. 
Their inspiration comes from the classical study of Reye congruence for K3 surfaces. 
Homological mirror symmetry and Fourier--Mukai partners are very important subjects but we refrain from going into details of this aspect of mirror symmetry in this article. 

Before closing this section, we note that there are various versions of mirror symmetry 
and the reader should not regard the formulations we have discuss above as being the final words defining mirror symmetry. 
Ultimately they may be only reflecting a few symptoms of mirror symmetry and an eventual mathematical definition may prove to be quite different from the one we study today.

%%%%%%%%%%%%%%%%%%%%%%%%%%%%%%%%%%%%%%%%%%%%%%%%%%%%%%%%%%%%%%%%%%%%%%%%%%%%%%%%%%%%%%%%%%%%%%%%%%%%%%%%%%%%%%
%%%%%%%%%%%%%%%%%%%%%%%%%%%%%%%%%%%%%%%%%%%%%%%%%%%%%%%%%%%%%%%%%%%%%%%%%%%%%%%%%%%%%%%%%%%%%%%%%%%%%%%%%%%%%%

\subsection{Dolgachev--Nikulin mirror symmetry for K3 surfaces}
It is important to keep concrete examples in mind, and so we next give a brief review of the Dolgachev--Nikulin mirror symmetry for K3 surfaces \cite{Dol}. 
Mirror symmetry for abelian surfaces is similar. 

Recall that a K3 surface $X$ is a simply-connected Calabi--Yau surface. We assume $X$ is projective throughout the article. 
The second cohomology $(H^2(X,\Z),\langle*,** \rangle)$ equipped with the cup product is isomorphic to the K3 lattice 
$$
\Lambda_{K3}=U^{\oplus 3} \oplus E_{8} (-1)^{\oplus 2},
$$ 
where $U$ is the hyperbolic lattice and $E_8$ is the root lattice of type $E_8$. 
The K3 lattice $\Lambda_{K3}$ is the unique even unimodular lattice of signature $(3, 19)$. 
It is also endowed with a weight-two Hodge structure
$$
H^2(X,\C)=H^{2,0}(X)\oplus H^{1,1}(X)\oplus H^{0,2}(X).
$$
Let $\Omega$ be a holomorphic volume on $X$. 
The space $H^{2,0}(X)\cong\C$ is generated by the class of $\Omega$, which we denote by the same $\Omega$. 

The Neron--Severi lattice $NS(X)$ and the transcendental lattice $T(X)$ of $X$ are primitive sublattices of $H^2(X,\Z)$ defined respectively by
\begin{align}
NS(X)&=\{ x\in H^2(X,\Z) \ | \ \langle x,\Omega \rangle=0 \},\notag \\
T(X)& =NS(X)^\bot_{H^2(X,\Z)}.\notag 
\end{align}
Here we extend the bilinear form $\langle *,** \rangle$ on $H^2(X,\Z)$ to that on $H^2(X,\C)$ $\C$-linearly.
The Neron--Severi lattice $NS(X)$ is of signature $(1,\rho(X)-1)$, where $\rho(X)$ is the Picard number of $X$. 
It is also identified with the lattice of algebraic $2$-cycles and isomorphic to the Picard lattice $\mathrm{Pic}(X)$, with isomorphism induced by the first Chern class map. 

The Mukai lattice of $X$ is defined to be
$$
H^*(X,\Z) = H^0(X,\Z) \oplus  H^2(X,\Z) \oplus H^4(X,\Z)
$$
endowed with the product
$$
(p,l,s) \cdot (p',l',s')=\langle l,l'\rangle - \langle p,s'\rangle -\langle s,p'\rangle, 
$$
As an abstract lattice, we have 
$$
H^*(X,\Z)  \cong U^{\oplus 4} \oplus E_{8} (-1)^{\oplus 2}
$$
The weight-two Hodge structure on $H^2(X,\Z)$ also extends in such a way that the $(1,1)$-part is given by 
$$
H^0(X,\C)\oplus H^{1,1}(X)\oplus H^4(X,\C).
$$ 
The classical Torelli theorem asserts that two K3 surfaces $X$ and $X'$ are isomorphic if and only if there exists a Hodge isometry 
$$
H^2(X,\Z) \cong H^2(X',\Z).
$$
This was first proved for the projective K3 surface by Pjateckii-Shapiro--Shafarevich \cite{PSSha} and later for non-projective ones by Burns--Rapoport \cite{BurRap}. 
The derived version of this, due to Mukai \cite{Muk2} and Orlov \cite{Orl}, asserts that two K3 surfaces $X$ and $X'$ are derived equivalent if and only if there exists a Hodge isometry
$$H^*(X,\Z) \cong H^*(X',\Z)$$ of the Mukai lattices. 

The notion of a lattice-polarization is introduced by Nikulin \cite{Nik1} in an attempt to extend the idea of the usual polarization by an ample line bundle. 
The following definition is due to Dolgachev \cite{Dol}.  
\begin{Def}
Let $X$ be a K3 surface and $\iota: M \rightarrow NS(X)$ a primitive lattice embedding, where $M$ is a non-degenerate even lattice of signature $(1,k)$ for $0 \le k \le 19$.  
A pair $(X, \iota)$ is an ample $M$-polarized K3 surface if $\iota(M)$ contains an ample class. 
We often omit the embedding $\iota$ in our notation. 
\end{Def}

\begin{Def}
Let $M$ as above. 
A family of K3 surfaces $\pi:\mathcal{X}\rightarrow B$ is ample $M$-polarized if there is a Zariski open subset $B^o\subset B$ with trivial local system $\mathbb{M} \subset R^2(\pi|_{B^o})_*\Z$ 
such that $\mathbb{M}_b \subset NS(X_b)$ induces an ample lattice polarization of $X_b$ for $b \in B^o$. 
\end{Def}

There exists a (coarse) moduli space of ample $M$-polarized K3 surfaces of dimension $20-\mathrm{rank}(M)$. 

For the K3 surfaces, Hodge number mirror symmetry seems trivial at first sight, since every K3 surface has the identical Hodge numbers. 
The complex moduli space and the K\"ahler moduli space are somewhat mixed, as they both live in $H^2(X,\C)$. 
The Dolgachev--Nikulin mirror symmetry for K3 surfaces can be formulated, not as an exchange of the Hodge numbers, 
but as an exchange of the algebraic lattices and the transcendental lattices as follows. 

Let $X$ be an ample $M$-polarized K3 surface for a lattice $M$ of signature $(1,\rho-1)$. 
We fix a primitive isotropic vector $f$ in the orthogonal complement $M^\perp_{\Lambda_{K3}}$ of $M$ inside $\Lambda_{K3}$. 
We introduce another important lattice by 
$$
N=(\Z f)^\perp_{M^\perp}/\Z f
$$
of signature $(1,19-\rho)$, which admits a natural primitive embedding into $\Lambda_{K3}$. 
Then the Dolgachev--Nikulin mirror of $X$ is defined to be an ample $N$-polarized K3 surface $Y$. 
Note that we have the following sublattices of the full rank in the K3 lattice $\Lambda_{K3}$ 
$$
M \oplus U \oplus N \subset \Lambda_{K3}. 
$$
This shows that a generic ample $M$-polarized K3 surface $X$ has 
$$
NS(X)\cong M, \ \ \ T(X)\cong U \oplus N
$$
while 
a generic ample $N$-polarized K3 surface $Y$ has 
$$
NS(Y)\cong N, \ \ \ T(Y)\cong U \oplus M.
$$ 
Therefore we observe that the algebraic and transcendental $2$-cycles of mirror K3 surfaces are interchanged up to the hyperbolic factor $U$. 

To be more precise we need to consider the Mukai lattice $H^*(X,\Z)$ of a K3 surface $X$.  
The algebraic cycles of the K3 surface $X$ form the following sublattice of the Mukai lattice 
$$
H^0(X,\Z)\oplus NS(X)\oplus H^4(X,\Z) \cong  NS(X) \oplus U. 
$$
Then the Dolgachev--Nikulin mirror symmetry exchanges the algebraic and transcendental lattices of the mirror K3 surfaces. 
This is the hallmark of mirror symmetry for K3 surfaces. 

Mirror symmetry involves a mirror map which identifies the complex moduli space of a K3 surface and the K\"ahler moduli space of the mirror K3 surfaces. 
The tube domain for the complexified K\"ahler moduli space of ample $M$-polarized K3 surfaces is given by
$$
V(M) =\{B+\sqrt{-1} \kappa \ | \ \kappa^2>0\}. 
$$
On the other hand, the period domain for the ample $N$-polarized K3 surfaces is given by 
$$
D(M\oplus U) =\{[\Omega]\in \PP((M\oplus U)\otimes \C) \ | \ \Omega^2=0, \ \langle \Omega,\bar{\Omega} \rangle>0 \}. 
$$
Note that the transcendental lattice of a generic ample $N$-polarized K3 surfaces is isomorphic $M\oplus U$. 
The mirror map for the Dolgachev--Nikulin mirror symmetry for K3 surfaces can be understood by the tube domain realization of the type IV symmetric domains: 
$$
D(M\oplus U) \cong V(M).
$$

%%%%%%%%%%%%%%%%%%%%%%%%%%%%%%%%%%%%%%%%%%%%%%%%%%%%%%%%%%%%%%%%%%%%%%%%%%%%%%%%%%%%%%%%%%%%%%%%%%%%%%%%%%%%%%
%%%%%%%%%%%%%%%%%%%%%%%%%%%%%%%%%%%%%%%%%%%%%%%%%%%%%%%%%%%%%%%%%%%%%%%%%%%%%%%%%%%%%%%%%%%%%%%%%%%%%%%%%%%%%%

\subsection{Mirror symmetry for varieties with effective $-K_X$}
It is classically known that there is a version of mirror symmetry for the Fano manifolds.  
We expect that such mirror symmetry should hold also for varieties with effective anti-canonical divisors \cite{Aur1} although we must drop expectation that mirror to be algebraic even near large volume limits.  

We consider a pair $(X,Z)$ consisting of a smooth variety $X$ and an effective anti-canonical divisor $Z\in |-K_X|$. 
Although a choice of an anti-canonical divisor is often implicit, we keep it as a part of data in this article. 
The basic idea is that the complement $X\setminus Z$ can be thought of as a log Calabi--Yau manifold 
as there exists a holomorphic volume form $\Omega$ on $X\setminus Z$ with poles of order one along $Z$. 
A different choice of $Z$ gives rise to a different log Calabi--Yau manifold $X\setminus Z$.  

\begin{Ex} \label{toric}
Let $X$ be a toric Fano $n$-fold and $Z$ the toric boundary, which is the complement of the dense torus $(\mathbb{C}^\times)^n \subset X$. 
Then $X\setminus Z=(\mathbb{C}^\times)^n$ carries a standard holomorphic volume form 
$$
\Omega=\wedge_{i=1}^{n}\sqrt{-1}d\log z_i=\wedge_{i=1}^{n}\sqrt{-1} \frac{dz_i}{z_i},
$$
where $(z_i)$ are the coordinates of $(\mathbb{C}^\times)^n$. 
The toric boundary is a canonical choice of an anti-canonical divisor for a toric Fano manifold and such a choice is often implicit. 
\end{Ex}

\begin{Def}
A Landau--Ginzburg model is a pair $(Y,W)$ of a K\"ahler manifold $Y$ and a holomorphic function $W:Y\rightarrow \mathbb{C}$, 
which is called a superpotential. 
\end{Def}

Just like the Calabi--Yau case, there are various formulations of mirror symmetry for varieties with effective anti-canonical divisors. 
It is probably Eguchi--Hori--Xiong that first noticed that there is an A-twist of the $\sigma$-model associated to a quasi-Fano manifold \cite{EHX}.
They showed that this theory is the same as the theory coming from a Landau--Ginzburg model. 
In this article, we focus on a mirror conjecture for a quasi-Fano manifold together with a smooth anti-canonical divisor (Katzarkov--Kontsevich--Pantev \cite{KKP}, Harder \cite{Har}). 
Here a quasi-Fano manifold $X$ is a smooth variety $X$ such that $|-K_X|$ contains a smooth Calabi--Yau member and $H^i(X,\mathcal{O}_X)=0$ for $i>0$.

\begin{Conj} \label{LG MS}
For a pair $(X,Z)$ of a quasi-Fano $n$-fold $X$ and a smooth anti-canonical divisor $Z \in |-K_X|$, 
there exists a Landau--Ginzburg model $(Y,W)$ such that
\begin{enumerate}
\item the superpotential $W:Y\rightarrow \C$ is proper, 
\item for a regular value $s\in \C$ of $W$, we have 
$$
\sum_{j}h^{n-i+j,j}(X)=h^i(Y,W^{-1}(s))
$$
\item the generic fibres of $W$ and the generic anti-canonical hypersurfaces in $X$ are mirror families of compact Calabi--Yau $(n-1)$-folds, 
\end{enumerate}
where $h^i(Y,W^{-1}(s))$ is the rank of the relative cohomology group $H^i(Y,W^{-1}(s))$. 
The pair $(Y,W)$ is called a mirror Landau--Ginzburg model of $(X,Z)$. 
\end{Conj}

The anti-canonical divisor $Z$ can be thought of as an obstruction for the quasi-Fano manifold $X$ to be a Calabi--Yau manifold 
and similarly the superpotential $W$ is an obstruction for the Floer homology of a Lagrangian torus in $X$ to be defined in the sense of Fukaya--Oh--Ohta--Ono \cite{FOOO, CO} as we will see in the next section. 
Mirror symmetry for Calabi--Yau manifolds can be thought of as a special case of this conjecture when there is no obstruction. 

\begin{Ex} \label{HVMirror}
Let $X=\mathbb{P}^1$ and $Z=\{0,\infty\}$ equipped with a toric K\"ahler form $\omega$.  
Then the mirror Landau--Ginzburg model of $(X,Z)$ is given by 
$$
(\mathbb{C}^\times, W(z)=z +\frac{q}{z}),
$$
where $q=\exp(-\int_{\mathbb{P}^1}\omega)$.  
One justification of this mirror duality is given by the ring isomorphism 
$$
\mathrm{QH}(\mathbb{P}^1)=\mathbb{C}[H]/(H^2-q) \cong \mathbb{C}[z^{\pm1}]/(z^2-q) =\mathrm{Jac}(W). 
$$
Here $\mathrm{QH}(\mathbb{P}^1)$ is the quantum cohomology ring of $\mathbb{P}^1$ and $\mathrm{Jac}(W)$ is the Jacobian ring of the superpotential $W$. 
This version of mirror symmetry was first shown by Batyrev for the toric Fano manifolds. 
\end{Ex}

In Conjecture \ref{LG MS} we assume that $X$ is quasi-Fano so that there exists smooth $Z \in |-K_X|$. 
However, it can be generalized to the case when $Z$ is mildly singular. %at the cost of

\begin{Conj}
The mirror of a pair $(X,Z)$ consisting of a smooth variety $X$ and its effective anti-canonical divisor $Z \in |K_X|$ is given by 
a Landau--Ginzburg model $(Y,W)$ consisting of a K\"ahler manifold $Y$ and a holomorphic function $W:Y \rightarrow \C$. 
A generic fiber of $W$ is mirror symmetric to $Z$ and some Hodge theoretic equality similar to Conjecture \ref{LG MS} (2) still holds. 
\end{Conj}

This conjecture in particular claims that the anti-canonical divisor $Z$ is smooth if and only if the superpotential $W$ is proper. 
A good example the reader can keep in mind is the following. 

\begin{Ex} \label{P^2}
For $X=\mathbb{P}^2$ with a toric K\"ahler form $\omega$, we define $Z_0$ to be the toric boundary.  
The mirror Landau--Ginzburg model of the pair $(X,Z_0)$ is given by 
$$
(Y_0=(\mathbb{C}^\times)^2, W_0(x,y)=x+y +\frac{q}{xy}),
$$
where $q=\exp(-\int_{H}\omega)$ for the line class $H$.    
A generic fiber of $W_0$ is an elliptic curve with 3 punctures. 

On the other hand, let  $Z_1$ be the union of a smooth conic and a line intersecting 2 points, and $Z_2$ a nodal cubic curve, and $Z_3$ a smooth cubic curve. 
Then the mirror Landau--Ginzburg model $(Y_i,W_i)$ of the pair $(X,Z_i)$ is a fiberwise partial compactification of $W_0$ 
such that a generic fiber of $W_i$ is an elliptic curve with $3-i$ punctures for $1 \le i \le 3$. 
Note that these are compatible with a version of mirror symmetry for punctured and singular Riemann surfaces. 
\end{Ex}

The geometry of the Fano/Landau--Ginzburg correspondence is expected to be as rich as geometry of Calabi--Yau mirror symmetry. 
For instance, the moduli theory of Landau--Ginzburg models should be mirror to the birational geometry of quasi-Fano varieties, and vice versa. 

As briefly mentioned in the previous section, homological mirror symmetry can be extended to other classes of varieties beyond Calabi--Yau manifolds.  
The most studied examples are probably the toric Fano manifolds. 
Because of asymmetry of Fano/Landau--Ginzburg mirror symmetry, there are two versions of homological mirror symmetry depending upon on which side we consider the A-model (or the B-model).   
The complex moduli space of a toric Fano manifold $\PP_\Delta$ is trivial 
while the K\"ahler moduli space is of dimension $b_2=\dim H^2(\PP_\Delta)$ because the K\"ahler structure depends on the value of the moment map when we realize $\PP_\Delta$ as a GIT quotient. 
The latter non-trivial moduli space corresponds to the complex moduli of the mirror Landau--Gibzburg model $(Y,W)$, that is, the coefficients of the superpotential $W$. 
One version of homological mirror symmetry is an equivalence of (enhanced) triangulated categories 
$$
\mathrm{D^bCoh}(X)\cong \mathrm{D^bFS}(W), 
$$
where $\mathrm{D^bFS}(W)$ denotes the Fukaya--Seidel category of the Lefshetz fibration induced by $W$. 
The other version is 
$$
\mathrm{D^bFuk}(X)\cong \mathrm{MF}(W), 
$$
where $\mathrm{MF}(W)$ denotes the category of matrix factorization of $W$ under suitable gradings. 
The former version seems easier than the latter in the sense that there is no moduli involved, 
but the latter version is also tractable because it is expected that both $\mathrm{D^bFuk}(X)$ and $\mathrm{MF}(W)$ are semi-simple (thus the direct sum of the category of the complex of vector spaces) for a generic moduli parameter. 

%%%%%%%%%%%%%%%%%%%%%%%%%%%%%%%%%%%%%%%%%%%%%%%%%%%%%%%%%%%%%%%%%%%%%%%%%%%%%%%%%%%%%%%%%%%%%
%%%%%%%%%%%%%%%%%%%%%%%%%%%%%%%%%%%%%%%%%%%%%%%%%%%%%%%%%%%%%%%%%%%%%%%%%%%%%%%%%%%%%%%%%%%%%

\section{Degenerations and DHT conjecture} \label{Deg and DHT}
In this section, we will briefly review degenerations of Calabi--Yau manifolds, focusing on Tyurin degenerations. 
Mirror symmetry is in general a conjecture about a Calabi--Yau manifold near a large complex structure limit, which is thought to be a maximal degeneration, in the complex moduli space.   
However, in this article, we are mainly concerned with another class of loci in the complex moduli space, where a Calabi--Yau manifold degenerates to a union of two quasi-Fano manifolds.

%%%%%%%%%%%%%%%%%%%%%%%%%%%%%%%%%%%%%%%%%%%%%%%%%%%%%%%%%%%%%%%%%%%%%%%%%%%%%%%%%%%%%%%%%%%%%

\subsection{Tyurin degenerations and Heegaard splittings}

\begin{Def}
A Tyurin degeneration is a degeneration $\pi:\mathcal{X}\rightarrow \D$ of Calabi--Yau manifolds over the unit disc $ \D=\{|z|<1\} \subset \C$, 
such that the total space $\mathcal{X}$ is smooth and the central fibre $\mathcal{X}_0=X_1 \cup_Z X_2$ is a union of two quasi-Fano manifolds $X_1$ and $X_2$ 
intersecting normally along a common anti-canonical divisor $Z \in |-K_{X_i}|$ for $i=1,2$. 
We always assume that the central fiber is the unique singular fiber of the family. 
When there is no confusion, we simply write $X\rightsquigarrow X_1 \cup_Z X_2$ for a generic fiber $X$. 
\end{Def}

In the complex moduli space the points where the Tyurin degenerations occur usually form $1$-dimensional loci while the LCSLs are always $0$-dimensional and thus isolated. 
The Tyurin degenerations form a nice class of degenerations of a Calabi--Yau manifold when we try to study mirror symmetry beyond LCSLs. 
We will see this for K3 surfaces in the next section. 

\begin{Ex} \label{deg quintic}
For $k=2,3,5$, let $\sigma_k \in H^0(\PP^4,\mathcal{O}(k))$ be a generic section and denote the corresponding hypersurface by $X_k=\sigma_k^{-1}(0)\subset \PP^4$.  
We define
$$
\mathcal{X}' = \{(x,t) \in \PP^4 \times \mathbb \D \ | \  (t\sigma_5 + \sigma_2\sigma_3)(x)=0 \}
$$
and the second projection $\pi': \mathcal{X}' \rightarrow \D$. 
It is a degeneration of a quintic Calabi--Yau 3-fold to the union $X_2 \cup X_3$ of a quadric hypersurface $X_2$ and a cubic hypersurface $X_3$ in $\PP^4$. 
However, the total space $\mathcal{X}'$ is singular along the smooth curve $C=(X_2 \cap X_3 \cap X_5) \times \{0 \}$. 
The singular loci are locally the product of a smooth curve and a $3$-dimensional ordinary double point. 
If we blow-up $\mathcal{X}'$ along $C$, then the exceptional locus is a $\PP^1 \times \PP^1$-bundle over $C$. 
We contract one of the rulings of the bundle to get $\mathcal{X}$ and construct a Tyurin degeneration $ \pi:\mathcal{X} \rightarrow \mathbb \DD$ such that the central fiber is $X_2 \cup \mathrm{Bl}_CX_3$. 
\end{Ex}

Conversely, we can ask when the union of two quasi-Fano manifolds are smoothable to a Calabi--Yau manifold. 
Kawamata and Namikawa in fact showed the log deformation theory of a normal crossings Calabi--Yau manifold is unobstructed, 
by using similar techniques for proving the Bogomolov--Tian--Todorov theorem. 
Note that this was first proven for K3 surfaces by Friedman \cite{Fri}. 
We state a version of the claim, which is slightly modified for our setting.  

\begin{Thm}[{Kawamata--Namikawa \cite[Theorem 4.2]{KN}}] \label{KN}
Let $X_1$ and $X_2$ be quasi-Fano manifolds and $Z \in |-K_{X_i}|$ a common smooth anti-canonical divisor for $i=1,2$. 
Assume that there exist ample class $L_i \in \mathrm{Pic}(X_i)$ which restrict to an ample divisor $L_1|_Z=L_2|_Z$ on $Z$.  
Then the union $X_1 \cup_Z X_2$ of $X_1$ and $X_2$ intersecting normally along $Z$ is smoothable to a Calabi--Yau manifold $X$ 
if and only if 
$$
N_{Z/X_1}\cong N_{Z/X_2}^{-1}. 
$$
This is called d-semistability. 
Moreover the resulting Calabi--Yau manifold $X$ is unique up to deformation. 
\end{Thm}

\begin{figure}[htbp]
 \begin{center} 
  \includegraphics[width=55mm]{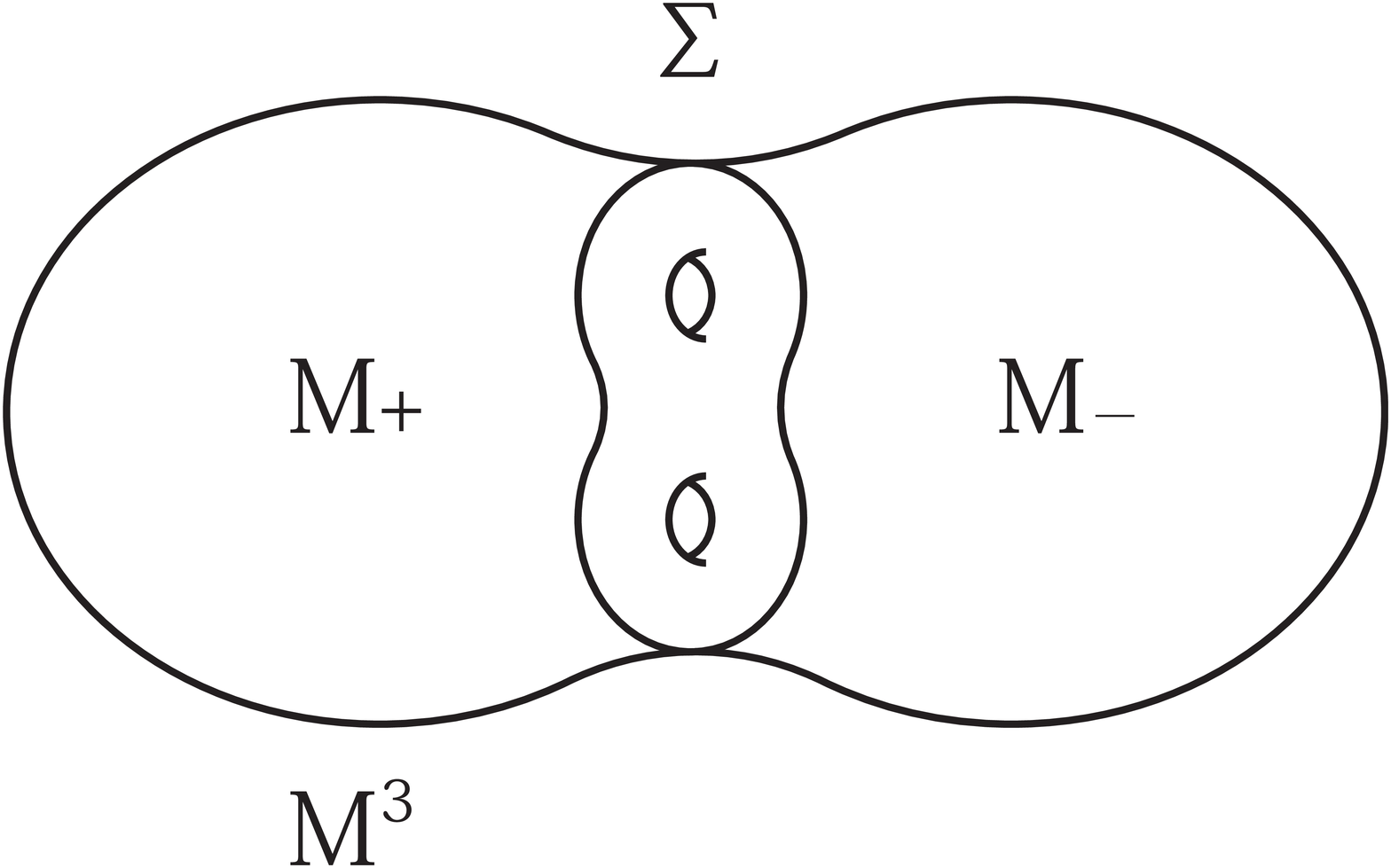}
 \end{center}%\vspace{-0.6zh}
  \caption{Heegaard splitting $M^3=M_+ \cup_\Sigma M_-$} 
\label{fig:Heegaard}
\end{figure}

A Tyurin degeneration of a Calabi--Yau 3-fold is thought to be a complex analogue of a Heegaard splitting $$M^3=M_+ \cup_\Sigma M_-$$ of a compact oriented real 3-fold without boundary $M^3$ 
by two handlebodies $M_{\pm}$ of genus $g$ with boundary a closed surface $\Sigma$, i.e. $\partial M_{\pm}=\pm \Sigma$ (Figure \ref{fig:Heegaard}).  
Indeed, a complex analogue of an oriented manifold without boundary is a Calabi--Yau manifold 
and that of an oriented manifold with boundary is a Fano manifold with an anti-canonical divisor. 
Note that a smooth anti-canonical divisor of a Fano manifold is a Calabi--Yau manifolds {\it because} the boundary of an oriented manifold is an oriented manifold without boundary. 

A Calabi--Yau manifold is called constructible if it is birational to a Calabi--Yau manifold which admits a Tyurin degeneration. 
The class of constructible Calabi--Yau manifolds is large enough to include all the complete intersection Calabi--Yau manifolds in the toric Fano manifolds. 
It is surprising that many rigid Calabi--Yau manifolds are in fact constructible. 
In his posthumous article \cite{Tyu}, Tyurin asked whether or not every Calabi--Yau 3-fold is constructible, just like every compact oriented real 3-fold admits a Heegaard splitting.  
Based on this analogy, he proposed to study geometry of a Calabi--Yau 3-fold by using that of quasi-Fano 3-folds when they are related by a Tyurin degeneration. 

Note that, although the motivation comes from 3-dimensional geometry, we do not restrict ourselves to 3-dimensions in this article.

%%%%%%%%%%%%%%%%%%%%%%%%%%%%%%%%%%%%%%%%%%%%%%%%%%%%%%%%%%%%%%%%%%%%%%%%%%%%%%%%%%%%%%%%%%%%%

\subsection{Degenerations of K3 surfaces} \label{deg K3}
To get a flavor of the Tyurin degenerations, let us take a close look at such degenerations of K3 surfaces. 
We begin with the Kulikov classification of the semi-stable degenerations of a K3 surface. 

\begin{Thm}[Kulikov \cite{Kul}, Persson--Pinkham \cite{PerPin}] \label{Kulikov}
Let $\pi: \mathcal{X} \rightarrow \D$ be a semi-stable degeneration of K3 surfaces with $K_{\mathcal{X}}=0$, such that all components of the central fiber $X_0=\pi^{-1}(0)$ are K\"ahler. 
Then one of the following occurs: 
\begin{enumerate}
\item (Type I) $X_0$ is a smooth K3 surface. 
\item (Type II) $X_0$ is a chain of elliptic ruled surfaces with rational surfaces at each end. 
The double curves of the end surfaces are anti-canonical divisors, and an elliptic ruled surface meets each nearby surface along an elliptic fiber. 
\item (Type III) $X_0$ consists of rational surfaces meeting along rational curves which form cycles in each component. 
The dual intersection comples of $X_0$ gives a triangulation of the 2-sphere $S^2$. 
\end{enumerate}
\end{Thm}
The cases in Theorem \ref{Kulikov} can also be distinguished by the action of monodromy on $H^2(X,\Z)$ of a general fibre $X$. 
Let $M$ denote the monodromy action on $H^2(X,\Z)$ around the origin $0 \in \D$.  
Then, after a base change if necessary, $N =\log M$ is nilpotent and has $N =0$ for Type I, $N \ne 0$ but $N^2 =0$ for Type II, and $N^2 \ne 0$ (and $N^3$=0 is automatic) for Type III.

We can also investigate the semi-stable degenerations of a K3 surface via the lattice theory. 
The Type II and III degenerations correspond respectively to the 0-dimensional and 1-dimensional cusps in the Baily--Borel compactification of the period domain of a K3 surface. 
Moreover there is a bijective correspondence between the 0-dimensional cusps and the primitive isotropic vectors in the transcendental lattice $T(S)$, up to automorphisms. 
Similarly for the 1-dimensional cusps and the 2-dimensional isotropic sublattices in $T(S)$, up to automorphisms. 

\begin{center}
 \begin{tabular}{|c|c|c|c|}  \hline
 & Type I & Type II & Type III \\ \hline
 $N=\log M$  & $N=0$ & $N\ne 0$, $N^2=0$ & $N^2 \ne 0$\\ \hline
moduli space  & smooth points & $0$-dim cusp & $1$-dim cusp \\ \hline
 $T(X)$  &  & $1$-dim isotropic & $2$-dim isotropic \\ \hline
 \end{tabular}
 \end{center}
 
A Tyurin degeneration is the simplest example of a Type II degeneration of a K3 surface when the central fiber $X_0$ has exactly two irreducible components. 
On the other hand, a Type III degeneration correspond to a LCSL of a K3 surface. 
Therefore, at least in the case of K3 surfaces, the Tyurin degenerations form a reasonable class of degenerations to study when we investigate mirror symmetry away from the LCSLs. 

Let us next examine the Tyurin degenerations (or more generally the Type II degenerations) from the viewpoint of the Dolgachev--Nikulin mirror symmetry. 
We begin with an ample $M$-polarized K3 surface $X$ for a lattice $$M \subset \Lambda_{K3} = U^{\oplus3}\oplus  E(-1)^{\oplus2}$$ of signature $(1,\rho-1)$. 
Once we fix a primitive isotropic vector $f \in M^\perp_{\Lambda_{K3}}$, 
we obtain the mirror of $X$ as an ample $N$-polarized K3 surface $Y$, where $$N=(\Z f)^\perp_{M^\perp}/\Z f$$ is a lattice of signature $(1,19-\rho)$.  

It is worth noting that this mirror correspondence depends upon the choice of an isotropic vector $f \in M^\perp_{\Lambda_{K3}}$. 
Fixing such an $f \in M^\perp_{\Lambda_{K3}}$ is equivalent to fixing a $0$-dimensional cusp $P$ 
(and thus a Type III degeneration) in the Baily--Borel compactification $\overline{\mathcal{D}(N\oplus U)}$ of the period domain $\mathcal{D}(N\oplus U)$ of the ample $M$-polarized K3 surfaces. 

In fact, based on the above investigation, we can show that there is a bijective correspondence between $1$-dimensional cusps in $\overline{\mathcal{D}(N\oplus U)}$ that pass through the cusp $P$, 
and the primitive isotropic vectors $e \in N$ up to isomorphism \cite{Dol,DHT}.  
Since a $1$-dimensional cusp in $\overline{\mathcal{D}(N\oplus U)}$ corresponds to a Type II degeneration of $X$, and an isotropic vector $e \in N \cong NS(Y)$, up to reflections, gives rise to an elliptic fibration on the mirror $Y$. 

The observation that {\it the Type II degenerations correspond to the elliptic fibrations on the mirror K3 surface} is compatible with the Doran--Harder--Thompson conjecture, which we will discuss below.

%%%%%%%%%%%%%%%%%%%%%%%%%%%%%%%%%%%%%%%%%%%%%%%%%%%%%%%%%%%%%%%%%%%%%%%%%%%%%%%%%%%%%%%%%%%%%
%%%%%%%%%%%%%%%%%%%%%%%%%%%%%%%%%%%%%%%%%%%%%%%%%%%%%%%%%%%%%%%%%%%%%%%%%%%%%%%%%%%%%%%%%%%%%

\subsection{DHT conjecture}  \label{DHT}
Mirror symmetry for Calabi--Yau manifolds and that for (quasi-)Fano manifolds have been studied for a long time, but somewhat independently. 
A natural question to ask is, 
\begin{center}{\it   
How are mirror symmetry for these manifolds related to each other?}.\end{center} 
As a matter of fact, when the author first leaned about mirror symmetry, he was puzzled by the asymmetry: 
a mirror of a Calabi--Yau manifold is a Calabi--Yau manifold while a mirror of a Fano manifold is a Landau--Ginzburg model. 
It was also hard to understand Landau--Ginzburg models geometrically. 

Motivated by works of Dolgachev \cite{Dol}, Tyurin \cite{Tyu}, and Auroux \cite{Aur2}, in String-Math 2016 held at the Tsinghua Sanya International Mathematics Forum, 
Doran--Harder--Thompson proposed the following geometric conjecture, which we call the DHT conjecture for short. 

\begin{Conj}[Doran--Harder--Thompson \cite{DHT, Kan1}]
Given a Tyurin degeneration of a Calabi--Yau manifold $X$ to the union $X_1 \cup_Z X_2$ of quasi-Fano manifolds intersecting along their common smooth anti-canonical divisor $Z$, 
then the mirror Landau--Ginzburg models $W_i:Y_i\rightarrow \C$ of $(X_i,Z)$ for $i=1,2$ can be glued together 
to be a Calabi--Yau manifold $Y$ equipped with a Calabi--Yau fibration $W:Y\rightarrow \PP^1$. 
Moreover, $Y$ is mirror symmetric to $X$. 
\end{Conj}

\begin{figure}[htbp]
 \begin{center} 
  \includegraphics[width=58mm]{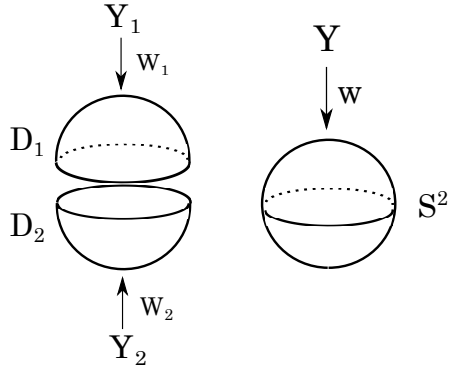}
 \end{center}
  \caption{Gluing Landau--Ginzburg models} 
\label{fig:gluing}
\end{figure}

The above gluing process can be understood as follows. 
We denote by $n$ the dimension of $X$ and by $Z_i^\vee$ a fiber of the superpotential $W_i$ mirror to a Calabi--Yau $(n-1)$-fold $Z$. 
\begin{enumerate}
\item Firstly, we assume that all the important information about the Landau--Ginzburg model $W_i:Y_i\rightarrow \C$ is contained in the critical locus of the superpotential $W_i$. 
Therefore, without much loss of information, we may replace it with a new Landau--Ginzburg model $W_i:Y_i\rightarrow \DD_i$ 
for a sufficiently large disc $\DD_i$ which contains all the critical values by shrinking $Y_i$ accordingly. 
\item Secondly, the Calabi--Yau manifolds $Z_1^\vee$ and $Z_2^\vee$ are both mirror symmetric to $Z$, and thus we expect that they are topologically identified. 
Note that two Calabi--Yau manifolds may be topologically different even if they share the same mirror manifold. 
There is no problem if $\dim X=1,2$ or $3$ for example.   
\item Thirdly, Theorem \ref{KN} implies that we have the d-semistability $N_{Z/X_1}\cong N_{Z/X_2}^{-1}$ because $X_1 \cup_Z X_2$ is smoothable to a Calabi--Yau manifold $X$. 
According to Kontsevich's homological mirror symmetry \cite{Kon}, we have an equivalence of triangulated categories
$$\mathrm{D^bCoh}(Z)\cong \mathrm{D^bFuk}(Z_i^\vee).$$ 
Then the monodromy symplectomorphism on $Z_i^\vee$ associated to the anti-clockwise loop $\partial \DD_i$ can be identified 
with the autoequivalence  $(-)\otimes \omega_{X_i}[n]|_Z$ on $ \mathrm{D^bCoh}(Z)$ (see \cite{Sei1, KKP} for details). 
By the adjunction formula, we have $$(-)\otimes \omega_{X_i}[n]|_Z\cong (-)\otimes N_{Z/X_i}^{-1}[n].$$ 
Therefore the d-semistability $N_{Z/X_1}\cong N_{Z/X_2}^{-1}$ implies, under mirror symmetry, that the monodromy action on $Z_1^\vee$
along the anti-clockwise loop $\partial \DD_1$ and that on $Z_2^\vee$ along the clockwise loop $-\partial \DD_2$ can be identified. 
\end{enumerate}
Therefore, assuming various mirror symmetry statements, 
we are able to glue the fibrations $W_i:Y_i\rightarrow \DD_i$ for $i=1,2$ along open neighborhoods of the boundaries $\partial \DD_1$ and $\partial \DD_2$ 
to construct a $C^\infty$-manifold $Y$ equipped with a fibration $W:Y\rightarrow S^2$ (Figure \ref{fig:gluing}).

Note that the smoothness of $Z$ implies the compactness of $Y$, which follows from the properness of the superpotentials $W_i$ for $i=1,2$ (c.f. Example \ref{P^2}). 
The highly non-trivial part of the conjecture is that there exists a Calabi--Yau structure on $Y$ and a complex structure on $S^2$ 
in such a way that $W:Y\rightarrow \PP^1$ is holomorphic and $Y$ is mirror symmetric to the Calabi--Yau manifold $X$. 

In \cite{DHT} the authors provide supporting evidence for this conjecture in various different settings, 
including the Batyrev--Borisov mirror symmetry and the Dolgachev--Nikulin mirror symmetry. 
For example, in the 3-fold case, topological mirror symmetry is proven to be equivalent to a version of Dolgachev--Nikulin mirror symmetry for K3 surfaces, provided that $Y$ admits a Calabi--Yau structure.  
Another important result is that, under reasonable assumptions, the resulting $C^\infty$-manifold $Y$ has the expected Euler number 
$$
\chi(Y)=(-1)^{\dim X}\chi(X).
$$   
Thus the conjecture is essentially proven at the topological level.  
However, the real difficulty of the conjecture lies in constructing $Y$ as a complex manifold, which should be mirror to the symplectic manifold $X$ (or vice versa). 

\begin{Rem}
It is worth mentioning that Calabi--Yau mirror symmetry and Fano/Landau--Ginzburg mirror symmetry look very different in physics.
In the former case the physical picture is an isomorphism of two conformal field theories while in the latter case it is an isomorphism of two massive theories. 
It is interesting to find physical aspects of the DHT conjecture. 
\end{Rem}

The DHT conjecture for the elliptic curves was recently proven by the author by using ideas of SYZ mirror symmetry. 
The proof is of interest in its own right and can be generalized to abelian varieties in a relatively straightforward way. 
We will include the proof in this article and also present an extension of the abelian surface case.

%%%%%%%%%%%%%%%%%%%%%%%%%%%%%%%%%%%%%%%%%%%%%%%%%%%%%%%%%%%%%%%%%%%%%%%%%%%%%%%%%%%%%%%%%%%%%
%%%%%%%%%%%%%%%%%%%%%%%%%%%%%%%%%%%%%%%%%%%%%%%%%%%%%%%%%%%%%%%%%%%%%%%%%%%%%%%%%%%%%%%%%%%%%

\subsection{Complex and K\"ahler degenerations}  \label{complex and Kahler deg}
In light of mirror duality, a complex degeneration of a Calabi--Yau manifold $X$ is mirror to a K\"ahler degeneration of a mirror Calabi--Yau manifold $Y$. 
For instance, Morrison proposed that geometric transitions are reversed under mirror symmetry \cite{Mor1,KanLau1}. 
Recall that a geometric transition is a birational contraction followed by a complex smoothing, or in the reverse way, applied to a K\"ahler manifold.   
Thus Morrison's conjecture claims that the mirror of a birational contraction is a certain complex degeneration. 

We now try to understand the DHT conjecture from a slightly different perspective, namely a K\"ahler degeneration.  
The author learned the following idea in discussion with Andrew Harder. 
If we think classically, every degeneration of a Calabi--Yau manifold should correspond to a contraction on the mirror side as explained above. 
Thus any Tyurin degeneration of a Calabi--Yau manifold $X$ should be mirror dual to a contraction of the mirror Calabi--Yau manifold $Y$.  
In our case the mirror contraction is so bad that it contracts a family of divisors, forming a {\it Calabi--Yau fibration} from the mirror Calabi--Yau manifold to $\PP^1$.  

It may be useful to give supporting evidence to this observation. 
As analyzed in \cite[Section 5.1]{DHT} we start with a Tyurin degeneration of a Calabi--Yau 3-fold $X$ and analyze the monodromy action $M$ of this degeneration on the cohomology group $H^3(X,\C)$. 
Here we assume that the Tyurin degeneration loci is connected to a LCSL. 
Picking up a LCSL is equivalent to choosing a mirror Calabi--Yau manifold $Y$, just like we fix an isotropic vector in the K3 mirror symmetry. 

Under mirror symmetry, the associated operator action $N=\log(M)$ on $H^3(X,\C)$ corresponds to 
the action of the cup product with $c_1(L)$ of a nef line bundle $L$ on $\oplus_{i=0}^3H^{\mathrm{i,i}}(Y,\C)$ of the mirror Calabi--Yau 3-fold $Y$: 
\begin{align}
N: & H^3(X,\C)  \longrightarrow H^3(X,\C), \notag  \\ 
c_1(L): & \oplus_{i=0}^3H^{\mathrm{i,i}}(Y,\C)  \longrightarrow \oplus_{i=0}^3H^{\mathrm{i,i}}(Y,\C). \notag
\end{align}
The assumption that we start with a Tyurin degeneration implies that $N \ne 0$ but $N^2=0$, and thus we have the mirror condition that $c_1(L) \ne 0$ but $c_1(L)^2=0$. 
Then it follows that $c_1(L)$ lies in the boundary of the K\"ahler cone of $Y$ and, by a result of Oguiso \cite{Ogu1}, 
some power of $L$ induces a map $Y \rightarrow \PP^1$ with fibers either abelian surfaces or K3 surfaces. 
This is nothing but the Calabi--Yau fibration mirror to the Tyurin degeneration in the DHT conjecture. 

\begin{Rem}
Around a LCSL, the logarithmic of a positive sum of the monodromy operators has a maximally unipotent monodromy 
i.e. $N^3 \ne 0$ ($N^4=0$ is automatic) and thus the mirror operator $c_1(L)$ is also maximally unipotent. 
This implies that $L$ lies inside the K\"ahler cone and gives a polarization of $Y$
\end{Rem}

In the K3 surface case, the Calabi--Yau fibration mirror to the Tyurin degeneration is  precisely the elliptic fibration we discussed in end of Section \ref{deg K3}.

The reader is warned that the DHT conjecture does not hold unless we impose a condition on the Tyurin degeneration of a Calabi--Yau manifold $X$.    
For example, if the complex moduli space of $X$ is $1$-dimensional, 
the K\"ahler moduli space of a mirror Calabi--Yau manifold $Y$ is also $1$-dimensional and thus $Y$ cannot have a fibration structure (unless it is 1-dimensional).  
The DHT conjecture should be modified so that the Tyurin degeneration occurs in a locus which contains a large complex structure limit. 
In such a case, the Calabi--Yau manifold $Y$ should be the mirror corresponding to the large complex structure limit. 
Otherwise what we could expect is that there exists a {\it homological mirror} $Y$ of $X$ 
equipped with a {\it non-commutative Calabi--Yau fibration} 
$$
W:\mathrm{D^bCoh}(\PP^1)\longrightarrow \mathrm{D^b}(Y)
$$ 
by Calabi--Yau categories $\mathrm{D^b}(Y) \otimes_{\mathrm{D^bCoh}(\PP^1)} \mathrm{D^bCoh}(p)$ for $p \in \PP^1$. 
This can be thought of as homological mirror to the Tyurin degeneration (see \cite[Section 6] {DHT} for more details). 

It is an interesting problem to understand this non-commutative fibration in the context of the K\"ahler moduli space.

%%%%%%%%%%%%%%%%%%%%%%%%%%%%%%%%%%%%%%%%%%%%%%%%%%%%%%%%%%%%%%%%%%%%%%%%%%%%%%%%%%%%%%%%%%%%%
%%%%%%%%%%%%%%%%%%%%%%%%%%%%%%%%%%%%%%%%%%%%%%%%%%%%%%%%%%%%%%%%%%%%%%%%%%%%%%%%%%%%%%%%%%%%%

\section{SYZ mirror symmetry} \label{SYZ}
The Strominger--Yau--Zaslow (SYZ) mirror symmetry conjecture \cite{SYZ} provides a foundational geometric understanding of mirror symmetry for Calabi--Yau manifolds. 
It claims that a mirror pair of Calabi--Yau manifolds should admit dual special Lagrangian torus fibrations. 

It is Hitchin \cite{Hit} who first observed that the base of the fibration, which is locally the moduli space of the special Lagrangian fibers \cite{McL}, carries two natural dual integral affine structures. 
These integral affine structures are essential in SYZ mirror symmetry and appear to be more fundamental than symplectic and complex geometry \cite{GroSie1}.   
One of the integral affine structures will play a vital role is our proof of the DHT conjecture.

%%%%%%%%%%%%%%%%%%%%%%%%%%%%%%%%%%%%%%%%%%%%%%%%%%%%%%%%%%%%%%%%%%%%%%%%%%%%%%%%%%%%%%%%%%%%%
%%%%%%%%%%%%%%%%%%%%%%%%%%%%%%%%%%%%%%%%%%%%%%%%%%%%%%%%%%%%%%%%%%%%%%%%%%%%%%%%%%%%%%%%%%%%%

\subsection{Special Lagrangian submanifolds} \label{sLag}

Let $X$ be a  symplectic $2n$-fold equipped a symplectic form $\omega$. 
A submanifold $L\subset X$ is called Lagrangian if $\omega|_L=0$ and $\dim L=n$.  
Given a smooth function $f:X\rightarrow \R$, we may deform a Lagrangian submanifold by the flow generated by the Hamiltonian vector field $V_f$, which is defined by $\omega(V_f,*)=df$ 
(more generally, the function $f$ can depend also on the time of the flow). 
Therefore the moduli space of the Lagrangian submanifolds is an $\infty$-dimensional space. 

Let us further assume that $X$ is a Calabi--Yau $n$-fold. 
The celebrated theorem of Yau \cite{Yau} asserts that $X$ admits the unique Ricci flat K\"ahler metric $g$ representing a specified K\"ahler class $\omega$.  
Then we can take a holomorphic volume form $\Omega$ by requiring 
$$
(-1)^{\frac{n(n-1)}{2}}\left(\frac{\sqrt{-1}}{2}\right)^{n}\Omega \wedge \overline{\Omega}=\frac{\omega^{n}}{n!}, 
$$
which determines $\Omega$ up to a phase $e^{\sqrt{-1} \theta} \in S^1$. 

\begin{Def}
A Lagrangian submanifold $L \subset X$ in a Calabi--Yau manifold equipped with a normalized holomorphic volume form $\Omega$ 
is called special of phase $\theta$ if $$\Im( e^{\sqrt{-1} \theta}\Omega)|_L=0$$ for some constant $\theta \in \R/2\pi\Z$. 
\end{Def}

Note that for a Lagrangian submanifold $L \subset X$, we have 
$$
\mathrm{vol}_L=e^{\sqrt{-1} \Theta}\Omega|_L. 
$$
for a function $\Theta:L\rightarrow \R/ 2\pi\Z$. 
Here $\mathrm{vol}_L$ denotes the Riemannian volume form on $L$ induced by the Ricci-flat metric associated to $\omega$. 
The function $\Theta$ is the Hamiltonian function for the mean curvature vector field $V_{MC}$, i.e. $\omega(V_{MC},*)=d\Theta$. 
This implies that the mean curvature flow preserves the Lagrangian condition and the function $\Theta$ satisfies the heat equation
$$
\frac{\partial \Theta}{\partial t}=\Delta \Theta.
$$ 
Heuristically the mean curvature flow moves a Lagrangian subamnifold $L$ to a union of possibly singular special Lagrangian submanifolds, 
on each irreducible component of which $\Theta$ is constant. 

The special Lagrangian submanifolds are calibrated manifolds and thus minimal submanifolds.  
The special conditions are much stronger than being merely Lagrangian, and in fact the moduli space of special Lagrangian submanifolds becomes finite dimensional. 
Moreover, McLean proved that the moduli space of special Lagrangian submanifolds was unobstructed, and identified its tangent space as follows. 

\begin{Thm}[McLean \cite{McL}] \label{McLean}
Let $X$ be a Calabi--Yau manifold and $L \subset X$ a compact special Lagrangian submanifold.  
The space of special Lagrangian deformations of $L \subset X$ is a manifold $B$, 
with tangent space at the point $[L] \in B$ corresponding to $M$ isomorphic to the space $\mathcal{H}^1(L,\R)$ of harmonic $1$-forms on $L$. 
\end{Thm}

A good way to think about special Lagrangian submanifolds is an analogy with the classical Hodge theory. 
We often consider Lagrangian submanifolds up to Hamiltonian isotopies and want to have good representatives. 
The situation is similar to considering the closed forms up to the exact forms, that is the de Rham cohomology, on a smooth compact manifold $M$. 
The Hodge theory asserts that given a Riemannian metric on $M$ we can choose the harmonic forms as good representatives of the equivalence classes: 
$$
\Ker(d)/\mathrm{Im}(d)|_{\Omega^k(M)}=H^k_{dR}(M,\R) \cong \mathcal{H}^k(L,\R)
$$
A folklore conjecture claims that we can choose special Lagrangian submanifolds as good representatives of Lagrangian submanifolds up to Hamiltonian isotopies. 
Here the additional data we put is a Calabi--Yau structure or $\Omega$. 
Theorem \ref{McLean} is the manifestation that this conjecture holds in the first order, but the situation is much more delicate in the global case (the representative may be singular and reducible).  
The uniqueness of special Lagrangian submanifolds in Hamiltonian deformation classes of Lagrangian submaniflds, under mild conditions was shown by Thomas and Yau \cite{ThoYau}.  
This subject is closely related to the stability conditions on the Fukaya category $\mathrm{D^bFuk}(X)$. 
A similar problem can be considered on the mirror side $\mathrm{D^bCoh}(Y)$ and such stability conditions were investigated by Douglas and Bridgeland \cite{Dou,Bri0,Bri1}.

A surjective map $\pi \colon X\rightarrow B$ is called a special Lagrangian fibration if a generic fiber is a connected special Lagrangian submanifold. 
By the Liouville--Arnold theorem, a fiber of a regular value is a torus if it is compact and connected.   

Again K3 surfaces provide good concrete examples of special Lagrnagian submanifolds. 
Let $X$ be a K3 surface equipped with a K\"ahler form $\omega$. 
By the celebrated Yau's theorem \cite{Yau}, there exists the unique Ricci-flat K\"{a}hler metric $g$ representing the class $[\omega]$.  
Then the holonomy group is $\mathrm{SU}(2)$ 
and the parallel transport defines complex structures $I,J,K$ satisfying the quaternion relations: $$I^{2}=J^{2}=K^{2}=IJK=-1$$ 
such that 
$$
S^{2}=\{aI+bJ+cK \in \End(TS) \ | \ a^{2}+b^{2}+c^{2}=1\}
$$ 
is the set of the possible complex structures for which $g$ is a K\"{a}hler metric.  
The period of $X$ in the complex structure $I$ is given by the normalized holomorphic volume form 
$$
\Omega_{I}(*,**)=g(J*,**)+\sqrt{-1}g(K*,**),
$$ 
and the compatible K\"{a}hler form is given by
$$
\omega_{I}(*,**)=g(I*,**).
$$ 
With respect to the complex structure $J,K$, the holomorphic volume forms and K\"{a}hler forms are respectively given by 
\begin{align}
\Omega_{J}=\omega_{I}+\sqrt{-1}\Re(\Omega_{I}) \ \ \ & \ \ \ \omega_{J}=\Im(\Omega_{I}), \notag \\
\Omega_{K}=\Im(\Omega_I)+\sqrt{-1}\omega_{I} \ \ \ & \ \ \ \omega_{K}=\Re(\Omega_{I}). \notag
\end{align}

The {\it hyperK\"{a}hler trick} asserts that a special Lagrangian $T^2$-fibration $\phi:X \rightarrow S^2$ with respect to the complex structure $I$ 
is the same as an elliptic fibration $\phi:X\rightarrow \PP^1$ with respect to the complex structure $K$. 
This is simply because a real smooth surface $S \subset X$ is holomorphic if and only if $\Omega|_S=0$ (Harvey--Lawson \cite{HL}). 
This suggests that the study of special Lagrangian torus fibration can be thought of as a vast generalization of the study of elliptic fibrations of K3 surfaces in high dimensions. 
It is an interesting and challenging problem to classify the singular fibers of special Lagrangian torus fibrations. 

It is important to keep in mind that a Calabi--Yau structure is controlled by the 3 real tensors $\Re(\Omega)$, $\Im(\Omega)$ and $\omega$. 
We should treat them on an equal footing, especially in the SYZ picture below.  
$$ 
\xymatrix{
\Re(\Omega) \ar@{<-->}[rd] & \ar@{-->}[l] \ar@{-->}[r] & \Im(\Omega) \ar@{<-->}[ld]\\
 & \omega & 
}
$$
In high dimensions, a holomorphic volume form $\Omega$ does not determine a complex structure of a Calabi--Yau manifold, 
but its local complex deformation is captured by the variations of $\Omega$.

%%%%%%%%%%%%%%%%%%%%%%%%%%%%%%%%%%%%%%%%%%%%%%%%%%%%%%%%%%%%%%%%%%%%%%%%%%%%%%%%%%%%%%%%%%%%%
%%%%%%%%%%%%%%%%%%%%%%%%%%%%%%%%%%%%%%%%%%%%%%%%%%%%%%%%%%%%%%%%%%%%%%%%%%%%%%%%%%%%%%%%%%%%%

\subsection{SYZ mirror symmetry for Calabi--Yau manifolds} \label{SYZ CY}
The celebrated SYZ mirror symmetry conjecture \cite{SYZ} asserts that, for a mirror pair of Calabi--Yau $n$-folds $X$ and $Y$, 
there exist special Lagrangian $T^n$-fibrations $\phi$ and $\phi^\vee$
$$ 
\xymatrix{
X \ar[rd]_\phi &  & Y \ar[ld]^{\phi^\vee}\\
 & B & 
}
$$
over the common base $B$, which are fiberwisely dual to each other away from singular fibers. 
This is motivated by the T-duality in string theory. 
The treatment of singular fibers constitutes the essential part of the conjecture where the quantum corrections come into the play.  
%It also suggests an intrinsic construction of the mirror $Y$ as the total space of the fiberwise dual the Lagrangian torus fibration on $X$. 

An importance of this conjecture lies in the fact that it leads to an intrinsic characterization of the mirror Calabi--Yau manifold $Y$;  
it is the moduli space of these special Lagrangian fibers $T^n \subset X$ decorated with flat $\mathrm{U}(1)$-connection. 
The SYZ conjecture not only provides a powerful tool to construct a mirror manifold $Y$ out of $X$ as a fiberwise dual, 
but also explains why mirror symmetry should hold via real Fourier--Mukai type transformations \cite{Leu}. 

It is worth noting that a mirror manifold $Y$ depends on the choice of a special Lagrangian fibration $\phi:X\rightarrow B$,  
and conjecturally this is equivalent to the choice of a large complex structure limit, 
where the Gromov--Hausdorff limit of the Calabi--Yau manifold $X$ with its Ricci-flat metric is topologically identified with the base $B$.  
This conjecture has been confirmed for certain K3 surfaces by Gross and Wilson \cite{GroWil2}, and later generalized to several other cases. 
We refer the reader to \cite[Section 6]{Gro3} for a nice comparison of complex degenerations and metric degeneration. 

Again we would like to take a look at the K3 surface $X$ case, where the choice of a LCSL is given by an isotropic vector $f \in T(X)$. 
This gives rise to a special Lagrangian torus fibration $\pi:X \rightarrow B \simeq S^2$, whose fiber class is $f$, up to reflections by root elements.   
If we choose a section $e$ of $\pi$, then they span a hyperbolic lattice $U \subset T(X)$, mirror to 
$$
H^0(Y,\Z)\oplus H^4(Y,\Z)\cong U
$$
in the Mukai lattice of the mirror K3 surface $Y$.  
This observation is compatible with homological mirror symmetry: 
the mirror of a special Lagrangian fiber is a skyscraper sheaf $\mathcal{O}_y$ of $Y$ and the mirror of a section is the structure sheaf $\mathcal{O}_Y$.

%%%%%%%%%%%%%%%%%%%%%%%%%%%%%%%%%%%%%%%%%%%%%%%%%%%%%%%%%%%%%%%%%%%%%%%%%%%%%%%%%%%%%%%%%%%%%

\begin{Def}
An integral affine structure on an $n$-dimensional real manifold $B$ is an atlas with transition functions in $\mathrm{GL}_n(\mathbb{Z})\ltimes \mathbb{R}^n$. 
Equivalently it is a collection of sections of the tangent bundle $TB$ which form a full rank fiberwise lattice $\Lambda \subset  TB$. 
\end{Def}

Let $\phi : X \rightarrow B$ be a special Lagrangian $T^n$-fibration of a Calabi--Yau $n$-fold $X$. 
We denote by $L_b$ the fiber of $\phi$ at $b \in B$. 
The complement $B^o \subset B$ of the discriminant locus carries two natural integral affine structures, which we call symplectic and complex. 
They are defined by $\omega$ and $\Im(\Omega)$ respectively as follows. 

Let $\{\gamma_i \} \subset H_1(L_b,\Z)$ and $ \{ \Gamma_i\} \subset H_{n-1}(L_b,\Z)$ be bases of the two homology groups.
 Then we can define $1$-forms $\{\alpha_i\} \subset \Gamma(B^o,T^*B^o)$ by defining their value on $v \in \Gamma(B^o,TB^o)$ as
$$
\alpha_i(v)=\int_{\gamma_i} \iota(v)\omega,
$$ 
where $\iota(v)$ denotes contraction by a lift of a tangent vector $v$ to a normal vector field of the corresponding special Lagrangian fiber. 
Similarly, we define $1$-forms $\{\beta_i\} \subset \Gamma(B^o,T^*B^o)$ by 
$$
\beta_i(v)=-\int_{\Gamma_i} \iota(v)\Im(\Omega).
$$
These forms are closed because $\omega$ and $\Omega$ are closed, and hence locally of the form $\alpha_i=dx_i$ and $\beta_i=d\check x_i$, with $\{x_i\}$ coordinates
defining the symplectic affine structure and $\{\check x_i\}$ defining the complex affine structure on the base $B^o$. 
Finally, the McLean metric $G$ on the base $B^o$, which should be defined by the Hessian of a potential function $K: B^o \rightarrow \RR$ in either two affine structures, is given by
$$
G(v_1,v_2)=-\int \iota(v_1)\omega\wedge \iota(v_2)\Im(\Omega).
$$
Here the integral is over the whole fiber. 
The symplectic and complex affine structures discussed above are Legendre dual to each other with respect to the potential $K$.

On the other hand, given an integral affine manifold $B$ of dimension $n$, we have smooth dual $T^n$-fibrations: 
$$ 
\xymatrix{
TB/\Lambda \ar[rd]_\phi &  & T^*B/\Lambda^* \ar[ld]^{\phi^\vee}\\
 & B & 
}
$$
Here $\Lambda=\mathbb{Z}\langle \frac{\partial}{\partial x_1},\dots,\frac{\partial}{\partial x_n}\rangle$ is a fiberwise lattice in the tangent bundle $TB$ 
generated by integral affine coordinates $\{x_i\}$ of $B$, and $\Lambda^*$ is the dual lattice in the cotangent bundle $T^*B$. 
In a natural way, $TB/\Lambda$ and $T^*B/\Lambda^*$ are complex and symplectic manifolds respectively. 

In order to make them (possibly non-compact) Calabi--Yau manifolds, we need a dual integral affine structure on $B$ 
so that the roles of $TB$ and $T^*B$ are swapped. 
More precisely, we need a potential function $K:B\rightarrow \mathbb{R}$ satisfying the real Monge-Amp\'ere equation 
$$
\det(\frac{\partial^2K}{\partial x_i\partial x_j} )=C
$$
for a constant $C \in \R$. 
Then the dual integral affine structure is given by the Legendre transformation of the original one. 
This is called semi-flat mirror symmetry and serves as a local model for SYZ mirror symmetry without quantum correction \cite{Leu}. 
In general it is a very hard problem to extend this picture when singular fibers are present. 

In Section \ref{DHT proof}, we will begin with a symplectic manifold $X$ and construct a complex manifold $Y$. 
So let us take a close look at this case. 
Given a special Lagrangian $T^n$-fibration $\phi:X\rightarrow B$ of a Calabi--Yau $n$-fold $X$, we endow $B^o$ with the symplectic integral affine structure. 
We may think of the semi-flat mirror $Y^o$ of $X^o=\phi^{-1}(B^o)$ as the space of pairs $(b,\nabla)$ 
where $b \in B^o$ and $\nabla$ is a flat $\U(1)$-connection on the trivial complex line bundle over $L_b$ up to gauge. 
There is a natural map $\phi^\vee:Y^o\rightarrow B^o$ given by forgetting the second coordinate. 
With the same notation as before, the complex structure of $Y^o$ is given by the following semi-flat complex coordinates
$$
z_i(b,\nabla)=\exp(-2\pi\int_{A_i}\omega)\Hol_{\nabla}(\gamma_i), 
$$ 
where $\mathrm{Hol}_{\nabla}(\gamma_i)$ denotes the holonomy of $\nabla$ along the path $\gamma_i$.  
Then we observe that the dual fibration $\phi^\vee$ is locally given by the tropicalization map 
$$
(z_i) \mapsto (-\frac{1}{2\pi}\log|z_i|)_i.
$$  
This is an analogue of a Lagrangian torus fibration in complex geometry, and gives fruiteful connections to non-archimedean geometry, tropical geometry and others.

%%%%%%%%%%%%%%%%%%%%%%%%%%%%%%%%%%%%%%%%%%%%%%%%%%%%%%%%%%%%%%%%%%%%%%%%%%%%%%%%%%%%%%%%%%%%%
%%%%%%%%%%%%%%%%%%%%%%%%%%%%%%%%%%%%%%%%%%%%%%%%%%%%%%%%%%%%%%%%%%%%%%%%%%%%%%%%%%%%%%%%%%%%%

\subsection{Superpotential via Fukaya category}  \label{SYZqFano}
This section will mostly follow the ideas of Auroux's fundamental article \cite{Aur1}. 
Let us consider a smooth variety $X$ of dimension $n$ with an effective anti-canonical divisor $Z \in |-K_X|$. 
Observing that the complement $X \setminus Z$ carries a holomorphic $n$-form with poles along $Z$, 
we think of $X \setminus Z$ as a log Calabi--Yau manifold, to which the above SYZ construction can be applied. 
Hence at least intuitively we can construct the SYZ mirror $Y^o$ of the complement $X \setminus Z$. 
However, the information about $Z$ is missing. 

The missing information is captured by the superpotential $W$ of a mirror Landau--Ginzburg model.  
In the frame work of SYZ mirror symmetry $W$ is obtained as the weighted count of holomorphic discs of Maslov index $\mu=2$ 
with boundary in a smooth fiber $L$ of a given special Lagrangian torus fibration $\phi:X \rightarrow B$ as follows. 
\begin{Def}
The superpotential $W$ is a function on the semi-flat mirror $Y^o$ given by  
$$
W(b,\nabla)=\sum_{\substack{\beta \in \pi_2(X,L_b) \\ \mu(\beta)=2}}n_\beta z_\beta(b,\nabla), 
$$
where $z_\beta$ is defined to be
$$
z_\beta(b,\nabla)=\exp(-2\pi\int_{\beta}\omega)\Hol_{\nabla}(\partial \beta)
$$
and $n_\beta$ denotes the one-point open Gromov--Witten invariant of class $\beta \in \pi_2(X,L)$ defined by the machinery of Fukaya--Oh--Ohta--Ono \cite{FOOO}. 
\end{Def}
In a good situation, the Maslov index is given by the formula $\mu(\beta)=2Z\cdot \beta$ (\cite[Lemma 3.1]{Aur1}). 
Therefore we observe the superpotential $W$ keeps track of the information about the anti-canonical divisor $Z$ that compactify $X\setminus Z$ to $X$.  
Moreover, a crucial observation is that the superpotential $W$ is locally a holomorphic function on $Y^o$. 

It is insightful to interpret the presence of holomorphic discs from the viewpoint of the Floer theory. 
Fukaya--Oh--Ohta--Ono \cite{FOOO} give an obstruction for the Lagrangian intersection Floer homology complex to be a genuine complex. 
Namely the double differential $\partial^2$ might not be zero and we instead consider a twisted version of the Fulkaya category. 
The usual Fukaya category is an $A_\infty$-category equipped with higher morphisms: 
$$
\mathfrak{m}_k: \Hom(L_0,L_1)\otimes \Hom(L_1,L_2)\otimes \dots \otimes \Hom(L_{k-1},L_{k}) \longrightarrow \Hom(L_1,L_k) 
$$
for $k \ge 1$. 
$\mathfrak{m_1}$ is the chain map $\partial$, $\mathfrak{m}_2$ is the multiplication associative only up to higher morphisms $\mathfrak{m}_k\ (k\ge2)$ 
which count holomorphic disks with their boundary on Lagrangian submanifolds $L_1,\dots,L_k$.  
For twisted Fukaya category we also have an obstruction term
$$
\mathfrak{m}_0: \C \rightarrow \Hom(L,L)
$$
which counts holomorphic disks with boundary on a single Lagrangian submanifold $L$. 
We usually considers not arbitrary Lagrangian submanifolds, but weakly unobstructed ones, 
namely those for which $\mathfrak{m}_0(1)$ is a scalar multiple of the strict unit $e_L$ of $\Hom(L,L)$. 
This happens, for example, when the minimal Maslov index of a holomorphic disc with boundary on $L$ is $2$ and Maslov index $2$ discs are regular. 
Then this $$\mathfrak{m}_0(1)=W \cdot e_L$$ gives nothing but the superpotential we defined above. 

\begin{Rem} \label{non toric}
It is also worth noting that the superpotential $W$ is in general not known to converge if $X$ is not a toric Fano manifold.  
Moreover, if a fiber Lagrangian submanifold $L$ bounds a holomorphic disc of class $\beta \in \pi_2(X,L)$ of Maslov index $0$, 
then the one-point open Gromov--Witten invariant $n_\beta$ depends on the fiber $L$ as well as the point $p \in L$ which the holomorphic discs are required to pass through. 
\end{Rem}

On the other hand, we often want to take a smooth anti-canonical divisor $Z$ instead of the toric boundary so that the mirror Landau--Ginzburg superpotential $W$ is proper. 
In this case there appears to be a discriminant locus in the interior of the base $B$ and we need quantum corrections in the above toric SYZ construction \cite{Aur1}.   
It is not clear whether or not the above version of the SYZ program works for the varieties with effective anti-canonical divisor. 

Finally we note that in Auroux's approach \cite{Aur1} the specialty condition is relative in the sense that it depends on the choice of an anti-canonical divisor. 
Also we do not know, for example whether or not there exists a special Lagrangian fibration on the compliment of a smooth cubic in $\PP^2$.

%%%%%%%%%%%%%%%%%%%%%%%%%%%%%%%%%%%%%%%%%%%%%%%%%%%%%%%%%%%%%%%%%%%%%%%%%%%%%%%%%%%%%%%%%%%%%
%%%%%%%%%%%%%%%%%%%%%%%%%%%%%%%%%%%%%%%%%%%%%%%%%%%%%%%%%%%%%%%%%%%%%%%%%%%%%%%%%%%%%%%%%%%%%

\subsection{SYZ mirror symmetry for toric Fano manifolds} \label{toric Fano}
 
To be more explicit and also to avoid the convergence issue (Remark \ref{non toric}), we shall focus on the toric Fano case. 
Namely, we consider a toric Fano $n$-fold $X$ equipped with a toric K\"ahler form $\omega$ and a meromorphic volume form $\Omega=\wedge_{i=1}^{n}\sqrt{-1}d\log z_i$, 
where $(z_i)_i$ are the standard coordinates of the open dense torus $(\C^\times)^n\subset X$. 
Let $Z\subset X$ be the the toric boundary (Example \ref{toric}). 

Then the toric moment map $\phi:X\rightarrow \R^n$ gives a smooth special Lagrangian $T^n$-fibration $\phi: X \setminus Z \rightarrow B^o$, 
where $B=\phi(X) \subset \R^n$ is the moment polytope and $B^o$ is its interior\footnote{By abuse of notation, we use the same $\phi$ for the restriction of $\phi$ to $X\setminus Z$. }. 
By construction of the semi-flat mirror, it is straightforward to check the following assertion. 

\begin{Prop} \label{semi-flat} 
We define the tropicalization map by 
$$
\mathrm{Trop}:(\C^\times)^n\longrightarrow \R^n, \ (z_i)_i\mapsto (-\frac{1}{2\pi}\log|z_i|)_i.
$$ 
Then the semi-flat mirror $Y^o$ of the complement $X\setminus Z \cong (\C^\times)^n$ is given by the polyannulus $\mathrm{Trop}^{-1}(B^o)$. 
Moreover the dual fibration $\phi^\vee$ is identified with the restriction $\phi^\vee=\mathrm{Trop}|_{Y^o}:Y^o\rightarrow B^o$. 
\end{Prop}
%%%%%%%%%%%%%%%%%%%%%%%%%%%%%%%%%%%%%%%%%%%%%%%%%%%%%%%%%%%%%%%%%%%%%%%%%%%%%%%%%%%%%%%%%%%%%

In the toric Fano case, we do not modify $Y^o$ further, so henceforth we simply write $Y=Y^o$. 
In general, there is a discriminant locus in the interior of $B$ and then the semi-flat mirror $Y^o$ needs quantum corrections by the wall-crossing formulae of the superpotential $W$.

Let us take a close look at the projective line $\PP^1$. 
We have a special Lagrangian $T^1$-fibration $\phi:\PP^1 \rightarrow B=[0,\Im(\tau)]$ given by the moment map, where $\tau=\sqrt{-1}\int_{\PP^1}\omega$. 
By Proposition \ref{semi-flat}, the mirror $Y$ of $\PP^1 \setminus \{0, \infty\} \cong \C^\times$ is given by the annulus 
$$
Y=A_{(q,1)}=\{ q< |z| < 1\} \subset \C,
$$
where $q=e^{2\pi \sqrt{-1}\tau}$. 
Each special Lagrangian $T^1$-fiber separates $\PP^1$ into two discs, one containing $0$ and the other containing $\infty$. 
The classes $\beta_1$ and $\beta_2$ representing these disc classes satisfy $\beta_1+\beta_2=[\PP^1]$, 
and hence the coordinates on $Y$ should satisfy $z_{\beta_1}z_{\beta_2}=q$. 
\begin{figure}[htbp]
 \begin{center} 
  \includegraphics[width=45mm]{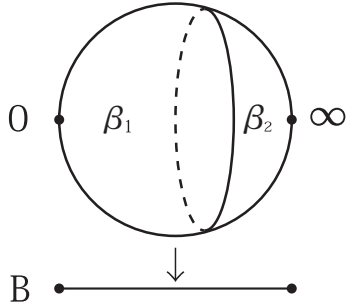}
 \end{center}
  \caption{$T^1$-fiber separates $\PP^1$ into two discs} 
\label{fig: P1}
\end{figure}

Moreover we can easily check that these are the only holomorphic discs of Maslov index $2$ and $n_{\beta_1}=n_{\beta_2}=1$. 
Using $z=z_{\beta_1}$ as a new coordinate on the mirror $Y=A_{(q,1)}$, 
we obtain the Landau--Ginzburg superpotential
$$
W(z)=z_{\beta_1}+z_{\beta_2}=z+\frac{q}{z}. 
$$
So far, we discuss only the real K\"ahler structure for simplicity, but we can easily complexify it in the above discussion.

%%%%%%%%%%%%%%%%%%%%%%%%%%%%%%%%%%%%%%%%%%%%%%%%%%%%%%%%%%%%%%%%%%%%%%%%%%%%%%%%%%%%%%%%%%%%%
%%%%%%%%%%%%%%%%%%%%%%%%%%%%%%%%%%%%%%%%%%%%%%%%%%%%%%%%%%%%%%%%%%%%%%%%%%%%%%%%%%%%%%%%%%%%%

\subsection{Renormalization} \label{renormalize}
The moment map for the toric $(S^1)^n$-action is defined only up to addition of a constant in the range $\mathrm{Lie}((S^1)^n)\cong \R^n$. 
In other words, the only intrinsic property of the base space $B$ is its affine structure and an affine embedding $B \subset \R^n$ is a choice. 
For example, we may take another moment map $\phi':\PP^1 \rightarrow B'=[\frac{-\Im(\tau)}{2},\frac{\Im(\tau)}{2}]$, and then the mirror Landau--Ginzburg model becomes
$$
W':Y'=A_{(q^{\frac{1}{2}},q^{-\frac{1}{2}})} \longrightarrow \C, \ \ \ z \mapsto z+\frac{q}{z},
$$
where 
$$
A_{(a,b)}=\{z \in \C \ | \ a < |z| < b\}
$$
for positive real numbers $a<b$. 
Note that we have a biholomorphism 
$$
A_{(q,1)} \cong A_{(q^{\frac{1}{2}},q^{-\frac{1}{2}})}, 
$$
which is induced by the translation of the underlying affine manifolds $B\cong B'$ in $\R^n$. 

Moreover, near the large volume limit, meaning $\int_{\PP^1}\omega  \gg 0$, we may identify $Y'$ with $\C^\times=A_{(0,\infty)}$. 
This renormalization procedure is also discussed by Hori and Vafa \cite{HV}. 
Auroux proposed a yet another renormalization as follows \cite[Section 4.2]{Aur1}. 
Let us consider a smooth variety $X$ with a nef anti-canonical divisor $Z\in |-K_X|$. 
Let $Y$ be the K\"ahler manifold which is the SYZ mirror of the complement $X\setminus Z$ as discussed above.  
We may enlarge $Y$ by using a renormalized K\"ahler form 
$$
\omega_k=\omega+k c_1(X)
$$
for large $k \gg 0$.  
Compared to the physical renormalization, this operation has the effect of not merely extending the domain $Y$ of the superpotential $W$, but also rescaling $W$ by the factor $e^{-k}$. 
However, by simultaneously rescaling the K\"ahler form and the superpotential $W$, we obtain a result consistent with the Hori--Vafa mirror.  

Anyway, we observe that the above construction reconstructs the Hori--Vafa mirror 
$$
(\C^\times,W(z)=z+\frac{q}{z}),
$$
after suitable renormalization (Example \ref{HVMirror}). 
We refer the reader to the influential articles \cite{HV, Aur1} for mode details of the renormalization process. 

\begin{Rem}
We will see that it is crucial in our proof of the DHT conjecture not to take the large volume limit but to keep track of the complex structures on the mirror annuli. 
In this way, we are able to naturally glue the Landau--Ginzburg models without the heuristic cutting process discussed in the DHT conjecture. 
\end{Rem}

%%%%%%%%%%%%%%%%%%%%%%%%%%%%%%%%%%%%%%%%%%%%%%%%%%%%%%%%%%%%%%%%%%%%%%%%%%%%%%%%%%%%%%%%%%%%%
%%%%%%%%%%%%%%%%%%%%%%%%%%%%%%%%%%%%%%%%%%%%%%%%%%%%%%%%%%%%%%%%%%%%%%%%%%%%%%%%%%%%%%%%%%%%%

\section{Degenerations and SYZ fibrations} \label{Deg and SYZ}

This section is logically not necessary, but we include it here to build a heuristic bridge between degenerations and SYZ fibrations of Calabi--Yau manifolds.    
It explains why mirror symmetry for Calabi--Yau manifolds only works near certain degeneration limits and why the SYZ conjecture needs to be viewed in a limiting sense.

%%%%%%%%%%%%%%%%%%%%%%%%%%%%%%%%%%%%%%%%%%%%%%%%%%%%%%%%%%%%%%%%%%%%%%%%%%%%%%%%%%%%%%%%%%%%%
%%%%%%%%%%%%%%%%%%%%%%%%%%%%%%%%%%%%%%%%%%%%%%%%%%%%%%%%%%%%%%%%%%%%%%%%%%%%%%%%%%%%%%%%%%%%%

\subsection{Approximating SYZ fibrations} \label{approximate SYZ}
Given a Calabi--Yau manifold, finding a special Lagrangian torus fibration is an important and currently unsolved problem in high dimensions, 
whereas there are some examples under relaxed conditions (see for example \cite{Gro2}).    
Among others, a well-known result is Gross and Wilson's work on Borcea--Voisin 3-folds \cite{Bor, Voi, GroWil2} and later it is generalized by the author and Hashimoto to Calabi--Yau 3-folds of type K \cite{HK1,HK2}. 
Even if we drop the speciality condition, finding a Lagrangian torus fibration of a given Calabi--Yau manifold is a very hard problem. 

\begin{figure}[htbp]
 \begin{center} 
  \includegraphics[width=55mm]{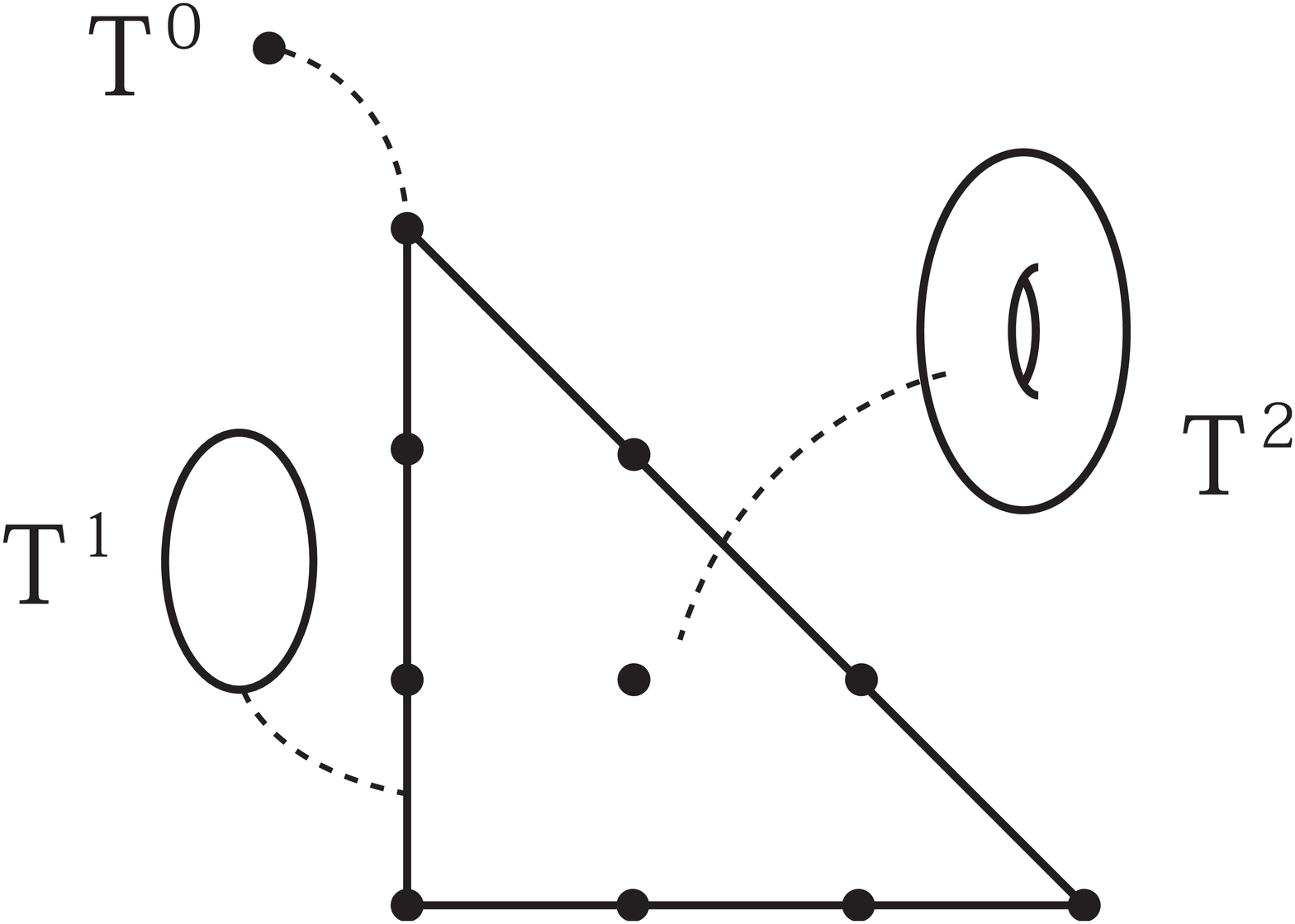}
 \end{center}
  \caption{Moment map for $\PP_\Delta=\PP^2$ polarized by $-K_{\PP^2}$} 
\label{fig:moment map}
\end{figure}

We would like to give some heuristic idea on constructing a Lagrangian torus fibration from a degeneration of a hypersurface Calabi--Yau manifold in a toric Fano manifold. 
The following idea is known among experts (c.f. \cite{LeuVaf,Rua}).  
Let $\Delta \subset M_\R$ be a reflexive lattice polytope (see Example \ref{Batyrev}).  
We consider the associated toric Fano $(n+1)$-fold $\PP_\Delta$ polarized by $-K_{\PP_\Delta}$, and choose a smooth Calabi--Yau manifold $X \in |-K_{\PP_\Delta}|$. 
The toric manifold $\PP_\Delta$ has the toric moment map $\pi:\PP_\Delta \rightarrow \Delta$, which is a Lagrangian $T^{n+1}$-fibration (Section \ref{toric Fano}).  
Recall that we have the fiber $\pi^{-1}(b) \cong T^{n+1-k}$ for a point $b$ in the strictly codimension $k$ boundary of $\Delta$ (Figure \ref{fig:moment map}).  
The moment polytope $\Delta$ encodes information about how tori $T^k$ of various dimensions are glued together to form the toric manifold $\PP_\Delta$. 
Moreover, $X_0=\pi^{-1}(\partial \Delta)$ is a union of toric manifolds, and the restriction $\pi|_{X_0}:X_0\rightarrow \partial \Delta$ is also a Lagrangian $T^{n}$-fibtation in an appropriate sense.  
 
Assume now that the Calabi--Yau manifold $X$ is close to $X_0=\pi^{-1}(\partial \Delta)$ in the linear system $|-K_{\PP_\Delta}|$. 
Let $r:X \rightarrow X_0$ be a continuous map induced by a degeneration.  
Then the composition 
$$
\phi=\pi|_{X_0}  \circ r:X\longrightarrow X_0 \longrightarrow \partial \Delta
$$
is thought to be an approximation of a Lagrangian $T^n$-fibration of $X$. 

In this approximate setting, we can understand an aspect of SYZ mirror symmetry as follows. 
Let $T_R^k$ denote a $k$-dimensional torus of radius $R>0$.  
We assume that we have the fiber 
$$
\phi^{-1}(b) \cong T_{1/\epsilon}^{n-k} \times T_\epsilon^k
$$
for small $\epsilon>0$ if a point $b$ lies in the strictly codimension $k$ strata of $\partial \Delta$. 
Following the idea of SYZ mirror symmetry, we dualize the fiber 
$$
(\phi^\vee)^{-1}(b) \cong T_\epsilon^{n-k} \times T_{1/\epsilon}^k
$$
by replacing $\epsilon$ by $\frac{1}{\epsilon}$ to obtain a mirror Calabi--Yau manifold $Y$.  
In toric geometry, this process is equivalent to taking the polar dual $\Delta^\vee \subset N_\R$ of the moment polytope $\Delta \subset M_\R$, 
replacing a $k$-dimensional face by an $(n-k)$-dimensional face. 
The base spaces $\partial \Delta, \partial \Delta^\vee$ are topologically the same, only the affine structures are different. 
These affine structures are related by a discrete version of Legendre transformation as we will see in Section \ref{toric deg}. 
Therefore we observe that the SYZ program provides us with a nice geometric interpretation of the Batyrev mirror construction (Example \ref{Batyrev}).

%%%%%%%%%%%%%%%%%%%%%%%%%%%%%%%%%%%%%%%%%%%%%%%%%%%%%%%%%%%%%%%%%%%%%%%%%%%%%%%%%%%%%%%%%%%%%
%%%%%%%%%%%%%%%%%%%%%%%%%%%%%%%%%%%%%%%%%%%%%%%%%%%%%%%%%%%%%%%%%%%%%%%%%%%%%%%%%%%%%%%%%%%%%

\subsection{Glimpse of toric degenerations} \label{toric deg}
The reader might notice that we do not really need the ambient toric variety $\PP_\Delta$ in the previous discussion. 
All we need is the property that the degeneration limit $X_0$ admits a Lagrangian torus fibration $\pi:X_0 \rightarrow B$. 
A naive way to construct such a Lagrangian torus fibration is simply gluing toric moment maps $\phi_\Delta:\PP_\Delta \rightarrow \Delta$.  
This leads us to the notion of toric degenerations. 

A degeneration of a Calabi--Yau manifold $\pi:\mathcal{X} \rightarrow \D$ is called toric 
if the central fiber $X_0$ is a union of toric varieties glued along toric strata and $\pi$ is locally toroidal (meaning a monomial in an affine toric chart).   
It is not necessarily semi-stable, and also in their program Gross and Siebert deal with more general situation where a generic fiber may have canonical singularities. 
Provided that irreducible components $\PP_\Delta$ of the central fiber $X_0$ have compatible toric polarization, 
we can construct a Lagrangian torus fibration $\pi:X_0 \rightarrow \cup_{\Delta}\Delta$, where we glue the lattice polytopes $\{\Delta\}_\Delta$ suitably. 
In the same manner as the previous section, for a Calabi--Yau manifold sufficiently close to $X_0$, 
we obtain an approximation of a Lagrangian torus fibration $\phi:X \rightarrow X_0 \rightarrow \cup_{\Delta}\Delta$. 
Note that the natural affine structure on the base $\cup_{\Delta}\Delta$ is the symplectic affine structure discussed in Section \ref{SYZ CY}. 

There are two important classes of examples. 
The first is toric degenerations of the toric varieties. 
Let $\PP_\Delta$ be the toric variety associated to a lattice polytope $\Delta \subset M_\R$, 
Given a regular subdivision $\Delta=\cup_{\Delta'}\Delta'$ of $\Delta$ which is the singular locus of some convex piecewise-linear function $\phi:\Delta \rightarrow \R$ with integral slopes, 
we define an unbounded polyhedron
$$
\widetilde{\Delta}=\{(m,h) \in M_\R \times \R \ | \ m \in M_\R, h \le \phi(m)\} \subset M_\R \times \R. 
$$
Then the second projection $\widetilde{\Delta} \rightarrow \R$ induces a toric morphism $\PP_{\widetilde{\Delta}}\rightarrow \C$, 
This is (after restriction to $\D$ for our definition) a toric degeneration such that a generic fiber is $\PP_\Delta$ 
and the central fiber is the union $\cup_{\Delta'}\PP_{\Delta'}$ of the toric varieties $\PP_{\Delta'}$ glued compatibly with the subdivision $\Delta=\cup_{\Delta'}\Delta'$.  
There is an obvious commutative diagram: 
$$ 
\xymatrix{
\PP_\Delta  \ar[d]_{\pi_\Delta} \ar@{~>}[r] & \cup_{\Delta'}\PP_{\Delta'} \ar[d]^{ \cup_{\Delta'}\pi_{\Delta'}}\\
\Delta  \ar@{~>}[r] &  \cup_{\Delta'}\Delta' 
}
$$
where $\phi_\Delta:\PP_{\Delta} \rightarrow \Delta$ is the toric moment map. 

The second is toric degenerations of the abelian varieties. 
Let $A$ be a $({\bf d})$-polarized abelian variety for ${\bf d} \in \Z^{n}_{\ge 1}$. 
The classical study of degenerations of the abelian varieties teaches us that any lattice subdivision of the $n$-dimensional torus $T^n=\R^n/{\bf d} \Z^n=\cup_{\Delta}\Delta$ gives rise to a toric degeneration 
such that a generic fiber is a $({\bf d})$-polarized abelian variety and the central fiber is the union $\cup_{\Delta}\PP_{\Delta}$ of toric varieties $\PP_{\Delta}$ (see \cite{OS, KanLau2}).    

\begin{figure}[htbp]
 \begin{center} 
  \includegraphics[width=45mm]{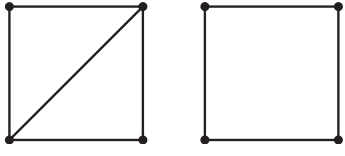}
 \end{center}
  \caption{Lattice subdivisions of $T^2=\R^2/\Z^2$} 
\label{fig: lattice subdiv}
\end{figure}

For instance, a principally polarized, that is $(1,1)$-polarized, abelian surface $A$ admits two classes of toric degenerations depending on whether we fully subdivide $T^2=\R^2/\Z^2$ or not (Figure \ref{fig: lattice subdiv}).  
In the former case we have the central fiber $\PP^2 \cup_{3\PP^1} \PP^2$ where we glue 3 $\PP^1$s according to the subdivision.  
In the latter case we obtain as the central fiber $\PP^1\times \PP^1/\sim$ where $\{0\}\times \PP^1 \sim \{\infty\}\times \PP^1$ and $\PP^1 \times\{0\} \sim \PP^1 \times \{\infty\}$. 

A $({\bf d})$-polarized abelian variety has a standard (special) Lagrangian torus fibration $\pi:A \rightarrow T^2$ and it is compatible with this class of toric degenerations 
$$ 
\xymatrix{
A  \ar[d]_{\pi} \ar@{~>}[r] & \cup_{\Delta}\PP_{\Delta} \ar[d]^{ \cup_{\Delta}\pi_{\Delta}}\\
T^2 \ar@{~>}[r] &  \cup_{\Delta}\Delta
}
$$

In the framework of SYZ mirror symmetry, an open set of the base of the SYZ fibration naturally carries an integral affine structure 
and the dual integral affine structure is obtained by the Legendre transformation of a potential function. 
In the setting of toric degenerations, the base space of the approximation of an SYZ fibration is the union of toric varieties $\cup_{\Delta}\Delta$. 
There is a discrete version of the usual Legendre transformation, 
given by the duality between a convex lattice polytope $\Delta \subset M_\R$ and a fan $\Sigma \subset N_\R$ together with a piecewise linear convex function $\psi$ on it: 
$$
\Delta \longleftrightarrow (\Sigma, \psi). 
$$
We readily see that $\psi(n)=\max\{\langle n,m\rangle | m \in \Delta\}$ is a discrete version of the Legendre transformation. 
This discrete Legendre duality applied to a union of lattice polytopes $\cup_{\Delta}\Delta$ also corresponds to the polar duality in toric geometry. 
More precisely, the discrete Legendre transformations is obtained by patching this example, and can be applied to affine manifolds with polyhedral decompositions. 

In the Gross--Siebert program \cite{GroSie0, GroSie1}, this duality is understood as a duality between the fan and cone pictures. 
They consider a polarized tropical manifold $(B,\mathcal{P},\phi)$, where $B$ is an integral affine manifold with singularities, 
$\mathcal{P}$ a nice polyhedral decomposition of $B$, and $\phi$ a strictly convex multivalued piecewise linear function on $B$. 
The program associates each polarized tropical manifold a toric degeneration of a Calabi--Yau manifold. 
They claim that Calabi--Yau manifolds associated to a polarized tropical manifold and its Legendre dual form a mirror pair.

%%%%%%%%%%%%%%%%%%%%%%%%%%%%%%%%%%%%%%%%%%%%%%%%%%%%%%%%%%%%%%%%%%%%%%%%%%%%%%%%%%%%%%%%%%%%%
%%%%%%%%%%%%%%%%%%%%%%%%%%%%%%%%%%%%%%%%%%%%%%%%%%%%%%%%%%%%%%%%%%%%%%%%%%%%%%%%%%%%%%%%%%%%%

\section{DHT conjecture via gluing} \label{DHT proof}

\subsection{DHT conjecture for elliptic curves} \label{Gluing}
We are in the position to confirm the DHT conjecture in the case of elliptic curves. 
In the language of SYZ mirror symmetry, we will {\it glue the Landau--Ginzburg models} by essentially gluing the affine base manifolds of special Lagrangian fibrations. 
The inspiration comes from the fact that a Tyurin degeneration is a complex analogue of a Heegaard splitting. 
Namely a Tyurin degeneration induces a Heegaard splitting of the base of the SYZ fibration: 
\begin{figure}[htbp]
 \begin{center} 
  \includegraphics[width=65mm]{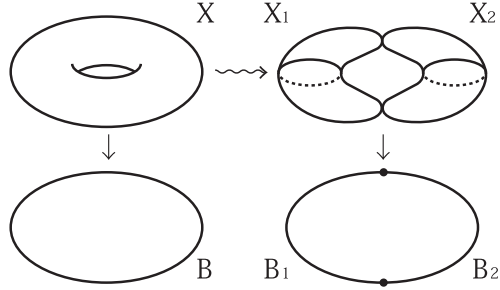}
 \end{center}
  \caption{SYZ fibration and Heegaard splitting} 
\label{fig: elliptic curve}
\end{figure}

Another key ingredient of the proof is to construct theta functions out of the Landau--Ginzburg superpotentials based on the geometry of the conjecture. 

Let us consider a Tyurin degeneration of an elliptic curve $\mathcal{X}\rightarrow \D$, 
where a generic fiber is a smooth elliptic curve and the central fiber $\mathcal{X}_0=X_1\cup_Z X_2$ is the union of two rational curves $X_1$ and $X_2$ glued at two points $Z$.  
We complexify the K\"ahler structure $\omega=B+\sqrt{-1} \kappa$ of $\mathcal{X}$ by introducing the B-field $B \in H^2(\mathcal{X},\R)/2 \pi H^2(\mathcal{X},\Z)$, 
and define 
$$
\tau_i=\int_{X_i}(B+\sqrt{-1}\kappa),\ \ \ q_i=e^{2 \pi \sqrt{-1} \tau_i}, \ \ \ i=1,2. 
$$ 
Then the complexified K\"ahler structure of a generic elliptic fiber $X$ of the family $\mathcal{X}\rightarrow \D$ is given by 
$$
\tau =\tau_1+\tau_2= \int_{\mathcal{X}_0} (B+\sqrt{-1} \kappa)
$$
so that $q=e^{2 \pi \sqrt{-1} \tau}=q_1q_2$.  

\begin{figure}[htbp]
 \begin{center} 
  \includegraphics[width=45mm]{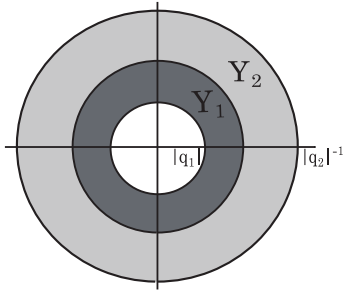}
 \end{center}
  \caption{$Y_1$ and $Y_2$} 
\label{fig:gluing annuli}
\end{figure}

Let us consider the moment maps $\phi_i:X_i \rightarrow B_i$ for $i=1,2$, 
where the base affine manifolds are $B_1=[0,\Im(\tau_1)]$ and $B_2=[-\Im(\tau_2),0]$. 
Then the mirror Landau--Ginzburgs of $(X_1,Z)$ and $(X_2,Z)$ are respectively given by 
\begin{align}
W_1:Y_1&=A_{(|q_1|,1)} \longrightarrow \C, \ \ z_1 \mapsto {z_1}+\frac{q_1}{z_1} \notag \\
W_2:Y_2&=A_{(1,|q_2|^{-1})} \longrightarrow \C, \ \ z_2 \mapsto z_2+\frac{q_2}{z_2}.  \notag
\end{align}
where $q_1z_2=z_1$.  
We observe that the boundary of the closure $\overline{Y_1 \cup Y_2} \subset \C^\times$ can be glued by the multiplication map $z \mapsto qz$ to form an elliptic curve $Y$, 
which is identified with the mirror elliptic curve $$\C^\times/q^{\Z}=\C^\times/(z \sim qz)$$ of $X$ (Figure \ref{fig:gluing annuli}).   
This construction corresponds to the gluing of the boundary of the affine manifold $B_1 \cup B_2=[-\Im(\tau_2),\Im(\tau_1)]$ 
by the shift $\Im(\tau)$ (twisted by the B-field upstairs).  

In order to confirm the DHT conjecture, we moreover want a map $\C^\times \rightarrow (\C \times \C)\setminus (0,0)$ 
which descends to a double covering $W:\C^\times/q^{\Z} \rightarrow \PP^1$ 
that locally looks like the superpotential $W_1$ (resp. $W_2$) over the upper (resp. lower) semisphere of the base $\PP^1$.   
Unfortunately, the naive analytic continuation of $(W_1,W_2)$ over $\C^\times$ does not work.  
The correct answer is given by considering all the Landau--Ginzburg models of the above sort, 
namely for $i \in \Z$ the Landau--Ginzburg models 
\begin{align}
W_{2i+1}:Y_{2i+1}&=A_{(|q^{-i}q_1|,|q^{-i}|)} \longrightarrow \C, \ \ z_{2i+1} \mapsto {z_{2i+1}}+\frac{q_{1}}{z_{2i+1}} \notag \\
W_{2i}:Y_{2i}&=A_{(|q^{1-i}|,|q^{1-i}q_2^{-1}|)} \longrightarrow \C, \ \ z_{2i} \mapsto z_{2i}+\frac{q_{2}}{z_{2i}}.  \notag
\end{align}
where the variable $z_i$ is defined inductively by $qz_{i+2}=z_i$. 
This means to eliminate the above arbitrary choice of a fundamental domain of the $\Z$-action on $\C^\times$ to construct the mirror elliptic curve $\C^\times/q^{\Z}$.   
A crucial observation is that if we consider all the even or odd superpotential $W_i$'s at once (in the sense of Remark \ref{product formula} below),  
they descend to the elliptic curve $\C^\times/q^{\Z}$ as sections of an ample line bundle as follows.  
Let us first consider the infinite product 
\begin{align}
W'_1(z) & =\prod_{i=1}^\infty (1+\frac{q_1}{z_{2i-1}^2})(1+\frac{z_{-2i+1}^2}{q_1}) \notag \\
&=\prod_{i=1}^\infty (1+q^{2i-1}(q_2z^2)^{-1})(1+q^{2i-1}q_2z^2) \notag \\
&= \prod_{k=1}^\infty \frac{1}{1-q^{2k}} \sum_{l \in \Z} q^{l^2}(q_2z^2)^l  \notag \\
&= \frac{e^{\frac{\pi \sqrt{-1} \tau}{6}}}{\eta(2\tau)}\vartheta_{\frac{1}{2},0}(2\zeta -\tau_1,2\tau), \notag
\end{align}
where we set $z=z_1=e^{2\pi \sqrt{-1} \zeta}$ and 
\begin{align}
\eta(\tau)=&e^{\frac{\pi \sqrt{-1} \tau}{12}}\prod_{m=1}^\infty(1-e^{2 \pi \sqrt{-1}\tau m}), \notag \\
\vartheta_{a,b}(\zeta,\tau)=&\sum_{n \in \Z} e^{\pi \sqrt{-1} (n+a)^2\tau} e^{2\pi \sqrt{-1} (n+a)(\zeta+b)} \notag
\end{align}
are the Dedekind eta function and theta function with characteristic $(a,b) \in \R^2$ respectively. 

\begin{Rem} \label{product formula}
We can regard $W'_1$ as the product of all the Landau--Ginzburg superpotential $W_i$'s for odd $i \in \Z$ because of the formula
$$
(z_j+\frac{q_k}{z_j})(z_{-j}+\frac{q_k}{z_{-j}})=q_kq^j (1+\frac{q_k}{z_j^2})(1+\frac{z_{-j}^2}{q_k}),
$$
for all $j \in \Z$ and $k=1,2$. 
\end{Rem}

In a similar manner, we next consider 
\begin{align}
W'_2(z) & =\prod_{i=1}^\infty (1+\frac{q_2}{z_{2i}^2})(1+\frac{z_{-2i+2}^2}{q_2}) \notag \\
&=\prod_{i=1}^\infty (1+q^{2i-1}\frac{q_1}{z^2})(1+q^{2i-1}\frac{z^2}{q_1}) \notag \\
&= \prod_{k=1}^\infty \frac{1}{1-q^{2k}} \sum_{l \in \Z} q^{l^2}(\frac{z^2}{q_1})^l  \notag \\
&= \frac{e^{\frac{\pi \sqrt{-1} \tau}{6}}}{\eta(2\tau)}\vartheta_{0,0}(2\zeta -\tau_1,2\tau), \notag
\end{align}
which can be thought of as the product of all the Landau--Ginzburg superpotential $W_i$'s for even $i \in \Z$. 
It is a classical fact that the theta functions with characteristics 
$$
\vartheta_{\frac{1}{2},0}(2\zeta -\tau_1,2\tau), \ \ \ \vartheta_{0,0}(2\zeta -\tau_1,2\tau)
$$
form a basis of the $(2)$-polarization of the elliptic curve $Y=\C^\times/q^{\Z}$. 
Therefore we obtain a double covering
$$
W:Y=\C^\times/q^{\Z} \longrightarrow \PP^1, \ z \mapsto [W'_1(z):W'_2(z)]. 
$$  
Observing that $W$ locally looks like the superpotential $W_i$ on each piece $Y_i$, 
we confirm that it is precisely the gluing of the two Landau--Ginzburg models argued in the DHT conjecture.  
This completes a proof of the conjecture in the case of elliptic curves.

It is interesting to observe that the product expressions of the theta functions are the manifestation of quantum corrections, 
which are encoded in the Landau--Ginzburg superpotentials, in SYZ mirror symmetry. 

The elliptic curves are somewhat special and our construction generalizes in a straightforward way to a degeneration of an elliptic curve to a nodal union of $n$ rational curves forming a cycle.   
The superpotential of each rational curve corresponds to a theta function with an appropriate characteristic, and they span a basis of the $(n)$-polarization of the mirror elliptic curve. 
We remark that the same result is obtained from different perspectives (c.f. \cite[Section 8.4]{Asp}, \cite[Section 4]{KanLau2}). 
This is due to the accidental fact that a Tyurin degeneration of an elliptic curve can be regarded as a degeneration at a large complex structure limit. 
A main difference shows up, for example, when we consider a type II degeneration of an abelian surface which is neither a maximal nor a toric degeneration (Section \ref{Gluing2}).   
However the essential mechanism of the DHT conjecture is already apparent in the case of elliptic curves: 
gluing the base affine manifolds and constructing theta functions from the Landau--Ginzburg superpotentials.

%%%%%%%%%%%%%%%%%%%%%%%%%%%%%%%%%%%%%%%%%%%%%%%%%%%%%%%%%%%%%%%%%%%%%%%%%%%%%%%%%%%%%%%%%%%%%
%%%%%%%%%%%%%%%%%%%%%%%%%%%%%%%%%%%%%%%%%%%%%%%%%%%%%%%%%%%%%%%%%%%%%%%%%%%%%%%%%%%%%%%%%%%%%
\subsection{DHT conjecture for abelian surfaces} \label{Gluing2}

We begin with a brief review of mirror symmetry for a class of abelian varieties.    
Let $E_\tau=\C/(\Z +\tau \Z)$ be an elliptic curve for $\tau \in \HH$. 
Then, for ${\bf d}=(d_i) \in \Z_{\ge 1}^n$, a ({\bf d})-polarized abelian variety is 
mirror symmetric to the self-product abelian $n$-fold $E_{d_1\tau}\times \dots \times E_{d_n\tau}$ (see for example \cite{KanLau2}).  

We shall discuss the surface case; 
for $(k,l)\in \Z_{\ge 1}^2$, a $(k,l)$-polarized abelian surface $X$ and a split abelian surface $Y$ of the form $E_{k\tau}\times E_{l\tau}$ are mirror symmetric. 
In this case we can understand the mirror duality in the same manner as the K3 surface case. 
A generic $(k,l)$-polarized abelian surface $X$ is an ample $\langle 2kl\rangle$-polarized abelian surface, 
where $\langle n \rangle$ denotes the lattice generated by $v$ with $v^2=n$. 
Then the transcendental lattice of $X$ is given by 
$$
T(X)\cong U^{\oplus2}\oplus \langle-2kl\rangle.
$$ 
On the other hand, a generic split abelian surface of the form $Y=E_{k\tau}\times E_{l\tau}$ is an ample $U\oplus\langle -2kl\rangle$-polarized abelian surface. 
Therefore $X$ and $Y$ are mirror symmetric in a sense similar to the Dolgachev--Nikulin mirror symmetry; 
namely the mirror duality comes from the following lattice splitting
$$
\langle 2kl\rangle \oplus U^{\oplus2}\oplus \langle-2kl\rangle \subset U^{\oplus3} \cong H^2(T^4,\Z). 
$$ 
Moreover, the author together with Fan and Yau recently investigated the Bridgeland stability conditions of $Y$ and identified its K\"ahler moduli space with the Siegel modular variety $\mathrm{Sp}(4,\Z) \backslash \mathfrak{H}_2$, 
which is nothing but the complex moduli space of $X$ \cite{FKY}.

Let us consider a Tyurin degeneration $\mathcal{X}\rightarrow \D$ of an $(k,l)$-polarized abelian surface $X$. 
It is a special case of a Type II degeneration and the degeneration limit is a union of two elliptic ruled surfaces $X_1 \cup_{Z} X_2$. (see \cite{HKW} for details of degenerations of polarized abelian surfaces). 
More precisely, there exists an elliptic curve $E$ and a degree $0$ line bundle $L \in \mathrm{Pic}^0(E)$ such that $$X_1\cong X_2\cong \PP(\mathcal{O}_E\oplus L)$$ 
and they are glued as follows. 
The $\infty$-section of $X_1$ is glued with the $0$-section of $X_2$. 
The $\infty$-section of $X_2$ is glued with the $0$-section of $X_1$ with translation by a gluing parameter $e \in E$. 
Thus we have $Z=2E$. 
Since $L$ is of degree $0$, the degeneration is topologically
$$
T^2 \times T^2  \rightsquigarrow (S^2 \cup_{\{0,\infty\}}S^2) \times T^2.
$$ 
By further deforming the degeneration we may assume that $X_i=\PP^1\times E$ for $i=1,2$ are glued along $\{0\}\times E$ and $\{\infty\}\times E$. 

\begin{Rem}
Strictly speaking, the above degeneration is not a Tyurin degeneration as the irreducible component $X_i$ is not a quasi-Fano manifold , 
but the DHT conjecture still makes sense.  
\end{Rem}

A $(k,l)$-polarized abelian variety $X$ admits two natural Lagrangian $T^2$-fibration. 
Namely, the Hodge form can be written as 
$$
\omega=kdx_1\wedge dx_3+ldx_2\wedge dx_4
$$
for a basis $\{dx_i\}_{i=1}^4$ of $H^1(X,\Z)$. 
We actually consider a complexification $\omega_\tau=\tau \omega$ for $\tau \in \mathbb{H}$. 
Then, for example, the projection in the $(x_1,x_3)$-direction gives rise a Lagrangian $T^2$-fibration $\pi:X \rightarrow T^2$. 
We assume that our Tyurin degeneration is compatible with this Lagrangian $T^2$-fibration.  
In other words, $\int_{\PP^1 \cup_{\{0,\infty\}} \PP^1}\omega_\tau=k\tau$ and $\int_{E}\omega_\tau=l\tau$ in the above notation. 

\begin{figure}[htbp]
 \begin{center} 
  \includegraphics[width=55mm]{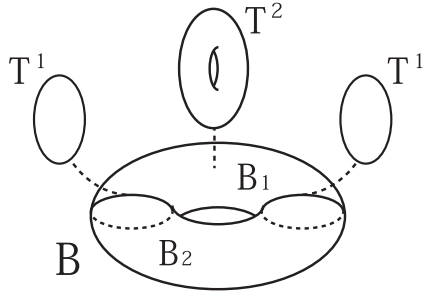}
 \end{center}
  \caption{Heegaard splitting induced by Tyurin degeneration} 
\label{fig:abel surface}
\end{figure}

\begin{Prop}
Let us consider the quasi-Fano manifold $\PP^1\times E$. 
Define $q=e^{2\pi \sqrt{-1} \int_{\PP1}\omega_\tau}$. 
Then the SYZ mirror of the pair $(\PP^1\times E,\{0,\infty\}\times E)$ is the Landau--Ginzburg model
$$
W:A_{(q,1)} \times E_{l \tau} \longrightarrow \C, \ z \mapsto z+\frac{q}{z}. 
$$
Note that the elliptic curve $E_{l \tau}$ is the mirror symmetric to $(E,\omega_\tau)$. 
\end{Prop}
This proposition is proven by a simple combination of SYZ mirror symmetry for a toric Fano manifold and that for an elliptic curve as we discussed in the the previous sections. 

Then it is not hard to generalize the argument in the previous section to this case. 
We can glue the two Landau--Ginzburg models of the above sort to construct the split abelian surface $Y=E_{k\tau} \times E_{l\tau}$ 
equipped with a Calabi--Yau fibration, whose fiber is given by the disjoint union of two $E_{l\tau}$: 
$$
W':Y=E_{k\tau} \times E_{l\tau} \stackrel{\mathrm{pr_1}}{\longrightarrow} E_{k\tau} \stackrel{2:1}{\longrightarrow} \PP^1. 
$$
Note that a generic fiber is $Z^\vee=2 E_{l\tau}$ mirror to $Z=2E$. 

The reader is also warned that $Y$ does not admit a Tyurin degeneration. 
It is also not hard to check the DHT conjecture for generic abelian surfaces, which are self-mirror symmetric. 
Note that an abelian surface $A$ splits if and only if the Neron--Severi lattice $NS(A)$ admits an embedding of the hyperbolic plane $U$. 
Therefore the DHT conjecture holds for a large classes of abelian surfaces. 

A key idea of our proof is to glue the two different affine manifolds $B_1$ and $B_2$ along the boundaries to obtain a compact affine manifold.  
This idea is not new and a similar construction (doubling) was already suggested by Auroux \cite{Aur2}. 
However, it is in general a very hard problem to glue together high dimensional affine manifolds along the boundaries.  
The difficulty is closely related to a choice of $Z \in |-K_X|$ of a quasi-Fano manifold $X$. 
More precisely, we need a Lagrangian fibration $\phi:X \rightarrow B$ which is compatible with $Z$, 
and then smoothness of $Z$ is likely to be proportional to that of the boundary $\partial B$. 
On the other hand, mirror symmetry for $(X,Z)$ tends to be harder as $Z$ gets less singular because we need to trade singularities of $\partial B$ with discriminant loci of the interior of $B$. 
This seems the main obstruction to generalizing our discussion in higher dimensions. 

Lastly, it is worth mentioning a recent work \cite{Lee} of Lee. 
Inspired by the DHT conjecture, he defined a dual pair of quasi-Fano manifolds which admit Calabi--Yau fibrations. 
He provided some evidence for his speculation that Calabi--Yau manfiolds $X$ and $Y$ which respectively admit Tyurin degenerations to $X_1 \cup X_2$ and $Y_1 \cup Y_2$ such that $X_i$ and $Y_i$ form a dual pair $i=1,2$ 
are mirror Calabi--Yau manifolds. 
From the viewpoint of the present article, his results can be understood as gluing of SYZ mirror symmetry of quasi-Fano manifolds which admit Calabi--Yau fibrations. 

%%%%%%%%%%%%%%%%%%%%%%%%%%%%%%%%%%%%%%%%%%%%%%%%%%%%%%%%%%%%%%%%%%%%%%%%%%%%%%%%%%%%%%%%%%%%%
%%%%%%%%%%%%%%%%%%%%%%%%%%%%%%%%%%%%%%%%%%%%%%%%%%%%%%%%%%%%%%%%%%%%%%%%%%%%%%%%%%%%%%%%%%%%%

\section{Splitting construction of Landau--Ginzburg models}  \label{split GL}
In this section we shall propose a construction of Landau--Ginzburg models. 
The inspiration comes from the DHT conjecture. 
The conjecture is originally about gluing but we will conversely consider splitting as described below.   
We note that large parts of that picture still remain conjectural. 

\subsection{Splitting conjecture}
Let us postulate that a Calabi--Yau manifold $X$ admits a Tyurin degeneration $X\rightsquigarrow X_1 \cup_Z X_2$ where $X_1$ and $X_2$ are quasi-Fano manifolds 
intersecting along a smooth anti-canonical divisor $Z$.  
Assume further that we know the mirror Calabi--Yau manifold $Y$ of $X$. 
As we discussed in Section \ref{complex and Kahler deg}, $Y$ should come equipped with a Calabi--Yau fibration $W:Y\rightarrow \PP^1$ corresponding to the Tyurin degeneration. 
In light of the DHT conjecture, we would like to propose the following. 

\begin{Conj} \label{split}
In the above setting, we are able to split the Calabi--Yau fibration $W$ into two fibrations 
$$
W_i=W|_{Y_i}:Y_i \longrightarrow \DD_i
$$
for $i=1,2$ along an embedded smooth loop $S^1 \subset \PP^1$, 
i.e. $\PP^1=\DD_1 \cup_{S^1}\DD_2$ and $Y_i=W^{-1}(\DD_i)$, in such a way that $(Y_i,W_i)$ is the mirror Landau--Ginzburg model of $(X_i,Z)$.  
\end{Conj}

The Conjecture \ref{split} involves geometry of Calabi--Yau manifolds. 
We need to analyze the moduli space of complex and K\"ahler moduli space to find the correspondence between Tyurin degenerations and Calabi--Yau fibrations. 
Moreover, we need to find a splitting whose monodoromy actions are compatible with the mirror Serre functors. 

In general it is a challenging problem to explicitly construct the mirror Landau--Ginzburg model of a given quasi-Fano manifold (or a variety with an effective anti-canonical divisor). 
Also, it is known that mirror symmetry for a quasi-Fano manifold tends to be harder as the anti-canonical divisor $Z$ gets less singular.  
The importance of the Conjecture \ref{split} lies in the fact it provides a completely new and explicit construction of mirror Landau--Ginzburg models which may be unreachable by the conventional approach.  
Moreover, it deals with the hardest case where the anti-canonical divisor $Z$ is smooth.

%%%%%%%%%%%%%%%%%%%%%%%%%%%%%%%%%%%%%%%%%%%%%%%%%%%%%%%%%%%%%%%%%%%%%%%%%%%%%%%%%%%%%%%%%%%%%

\subsection{Mirror symmetry for del Pezzo surfaces}
Let us begin with a brief review of mirror symmetry for del Pezzo surfaces. 
A (weak) del Pezzo surface $S$ is a smooth surface with ample (nef and big) anti-canonical divisor $K_S$. 
Such a surface has a degree $d$ defined as the self-intersection $K_S^2$. 
The possible degrees run between $1 \le d \le 9$. 
If $d \ne 8$, a (weak) del Pezzo surface of degree $d=9-k$ is the blow-up $\mathrm{Bl}_k\PP^2$ of $\PP^2$ at (almost) generic $k$ points. 
But if $d=8$, there are two possibilities: $\mathrm{Bl}_1\PP^2$ and $\mathbb{P}^1 \times \mathbb{P}^1$. 
Conversely the blow-up $\mathrm{Bl}_k\PP^2$ of $\PP^2$ at (almost) generic $k$ points is a (weak) del Pezzo surface for $0 \le k \le 8$. 
We refer the reader to \cite{Dem} for more details. 

In the seminal article \cite{AKO}, Auroux--Katzarkov--Orlov studied a version of homological mirror symmetry for the del Pezzo surfaces. 
Although their work deals with the del Pezzo surfaces, the mirror construction can be applied to the rational elliptic surface $\mathrm{Bl}_9\PP^2$ without modification. 
Let $\mathrm{Bl}_k\PP^2$ ($0 \le k \le 9$) be either a del Pezzo surface or an elliptic rational surface. 
Then they showed that the mirror Landau--Ginzburg model of $\mathrm{Bl}_k\PP^2$, equipped with a smooth $Z\in |-K_{\mathrm{Bl}_k\PP^2}|$, is a certain elliptic fibration 
$$
W_{k}:M_k\longrightarrow \C
$$
with $3+k$ singular fibers of Kodaira type $I_1$, 
which admits a compactification to a rational elliptic surface $\overline{W}_{k}:\overline{M}_k\rightarrow \PP^1$. 

It is worth mentioning a recent relevant work of Doran--Thompson \cite{DorTho} on mirror symmetry for lattice polarized (weak) del Pezzo surfaces. 
By introducing lattice polarizations, they showed that mirror symmetry for del Pezzo surfaces gains a significantly rich structure compared to the non-polarized case \cite{AKO}. 
Different choices of lattice polarization lead to different configurations of singular fibers (which may be deformed to the above fibration $W_k$). 

%%%%%%%%%%%%%%%%%%%%%%%%%%%%%%%%%%%%%%%%%%%%%%%%%%%%%%%%%%%%%%%%%%%%%%%%%%%%%%%%%%%%%%%%%%%%%

\subsection{Mirror symmetry for rational surfaces via K3 surfaces}
Probably a Tyurin degeneration of a K3 surface is the first non-trivial case to check the Conjecture \ref{split}. 
It is simple enough to be tractable but complex enough to display some essential features of the conjecture\footnote{
Everything in this section should work with $\mathbb{P}^1 \times \mathbb{P}^1$ instead of $\mathrm{Bl}_1\PP^2$ with slight modifications.  
For instance, a similar computation to Example \ref{deg quartic} can be carried out for a degeneration of a quartic K3 surface to a union of two quadric surfaces, which are isomorphic to $\PP^1\times \PP^1$.}.  

Let us first consider an ample $E_8^{\oplus 2}\oplus U$-polarized K3 surface $X$. 
Let $X_1 \cong X_2$ be rational elliptic surfaces $\mathrm{Bl}_9\PP^2$ obtained by blowing-up $\PP^2$ at the $9$ base points of a pencil of cubic curves. 
Then $X$ admits a natural Tyurin degeneration $X\rightsquigarrow X_1 \cup_Z X_2$ where $Z$ is a generic fiber of the elliptic fibration. 
This degeneration corresponds to the unique isotropic sublattice of rank $2$, up to isomorphism, in the transcendental lattice $T(X)\cong U^{\oplus 2}$ (Section \ref{deg K3}).  

Then the Dolgachv--Nikulin mirror K3 surface $Y$ of $X$ is an ample $U$-polarized K3 surface. 
It is nothing but a generic elliptic K3 surface with a section and thus comes equipped with an elliptic fibration $W:Y\rightarrow \PP^1$ mirror to the Tyurin degeneration.  
The fibration $W$ has 24 Kodaira $I_1$ fibers. 

Since the normal bundle $N_{Z/X_i}$ is trivial, the mirror monodromy symplectomorphism should be trivial. 
This implies that $W$ must split into two Landau--Ginzburg models $W_i:Y_i\rightarrow \DD_i$ for $i=1,2$ so that each $W_i$ has 12 singular fibers of Kodaira type $I_1$ 
(because of the famous 12 property\footnote{If each matrix $M_i$ is conjugate to $\begin{bmatrix} 1 & 1\\ 0 & 1\\ \end{bmatrix}$ and $\prod_{i=1}^kM_i=id$, then $12|k$.}).  
Then the Conjecture \ref{split} claims that $W_i:Y_i\rightarrow \DD_i$ is the mirror Landau--Ginzburg model of the rational elliptic surface $X_i$. 
We observe that our mirror Landau--Ginzburg model is compatible with the one constructed by Auroux--Katzarkov--Orlov \cite{AKO}.

\begin{Rem}
The difference between our mirror Landau--Ginzburg models and Auroux--Katzarkov--Orlov's will be explained by the renormalization procedure discussed in Section \ref{renormalize}. 
\end{Rem} 

Now that we consider a more general situation. 
We can construct a Tyurin degeneration $X\rightsquigarrow X_1 \cup_Z X_2$ of a K3 surface $X$ such that 
$X_1=\mathrm{Bl}_{k_1}\PP^2$ and $X_2=\mathrm{Bl}_{k_2}\PP^2$ for $k_1+k_2=18$ and $k_1,k_2\ge 3$, where the blow-ups occur at (almost) generic points. 
A simple example is given by a surface version of Example \ref{deg quintic}.  

\begin{Ex} \label{deg quartic}
For $k=1,3,4$, let $\sigma_k \in H^0(\PP^4,\mathcal{O}(k))$ be a generic section and we denote the corresponding hypersurface by $S_k=\sigma_k^{-1}(0)\subset \PP^3$.  
We define
$$
\mathcal{X}' = \{(x,t) \in \PP^3 \times \mathbb \D \ | \  (t\sigma_4 + \sigma_1\sigma_3)(x)=0 \}
$$
and the second projection $\pi': \mathcal{X}' \rightarrow \D$. 
It is a degeneration of a quartic K3 surface $S_4$ to the union $S_1 \cup S_3$ of a hyperplane $S_1\cong \PP^2$ and a cubic hypersurface $S_3 \subset \PP^3$. 
The total space $\mathcal{X}'$ has $12$ ordinary double points at $P=(S_1 \cap S_3 \cap S_4) \times \{0 \}$. 
We blow-up $\mathcal{X}'$ at $P$, then the exceptional locus is a $\PP^1 \times \PP^1$.  
We can contract one of the rulings to get $\mathcal{X}$ and construct a Tyurin degeneration $ \pi:\mathcal{X} \rightarrow \DD$ 
in such a way that the central fiber is $X_1 \cup_Z X_2$, where $X_1=\mathrm{Bl}_{k_1}\PP^2$ and $X_2=\mathrm{Bl}_{k_2}\PP^2$ for any $k_1+k_2=18$ and $k_2 \ge 6$. 
\end{Ex}

Let $W:Y\rightarrow \PP^1$ be the elliptic fibration mirror to the Tyurin degeneration $X\rightsquigarrow X_1 \cup_Z X_2$. 
We may assume that $X_1=\mathrm{Bl}_k\PP^2$ for $0 \le k \le 9$. 
Then the rational surface $X_2$ is not necessarily a weak del Pezzo surface nor an elliptic rational surface. 
Let $W_{k}:M_k \rightarrow \DD_1$ be the mirror mirror Landau--Ginzburg model of $X_1$ discussed by Auroux--Katzarkov--Orlov \cite{AKO} and Doran--Thompson \cite{DorTho}. 
In light of Conjecture \ref{split}, we claim that the mirror Landau--Ginzburg model of $X_2$ is given by an elliptic fibration 
$$
W^c_{k}:M^c_k \longrightarrow \DD_2
$$
obtained as the complement of the $W_{k}:M_k \rightarrow \DD_1$ in the elliptic fibration $W:Y\rightarrow \PP^1=\DD_1 \cup_{S^1}\DD_2$. 
Here we set $W_k=W|_{M_k=W^{-1}(\DD_1)}$ and $W^c_{k}=W|_{M^c_k=W^{-1}(\DD_2)}$. 

The singular fibers of $W^c_{k}$ may not be of the Kodaira type $I_1$.    
Also, the above fibration $W_{k}:M_k\rightarrow \C$ may not admit a compactification to a rational elliptic surface, in contrast to the construction of Auroux--Katzarkov--Orlov.

%%%%%%%%%%%%%%%%%%%%%%%%%%%%%%%%%%%%%%%%%%%%%%%%%%%%%%%%%%%%%%%%%%%%%%%%%%%%%%%%%%%%%%%%%%%%%

\subsection{Discussion}

Regarding to Conjecture \ref{split}, there appear two interesting and important problems. 
The first problem is to construct by our method high dimensional Landau--Ginzburg models which are not be tractable by the conventional approach. 
A very little is known about high dimensional Landau--Ginzburg models other than toric ones, especially those with proper superpotentials.   
The main difficulty of the problem lies in our poor understanding of an explicit correspondence between complex and K\"ahler degenerations in high dimensions.  
Probably the first case to be checked in 3-dimensions is Schoen's Calabi--Yau 3-fold \cite{Sch, HSS}, 
for which we may be able to carry out a calculation similar to an $E_8^{\oplus 2}\oplus U$-polarized K3 surface in the previous section. 
The new perspective of Conjecture  \ref{split} hopefully yields a new powerful way to construct Landau--Ginzburg models outside of the toric setting in high dimensions. 

The second problem is to consider a non-commutative version of Conjecture \ref{split}. 
Recall from the end of Section \ref{complex and Kahler deg} that, 
if a Tyurin degeneration $X\rightsquigarrow X_1 \cup_Z X_2$ of a Calabi--Yau manifold $X$ occurs in a locus which does not contain a LCSL, 
all we could expect is that there exists a {\it homological mirror} $Y$ of $X$ 
equipped with a {\it non-commutative Calabi--Yau fibration} 
$$
W:\mathrm{D^bCoh}(\PP^1)\longrightarrow \mathrm{D^b}(Y)
$$ 
by Calabi--Yau categories $\mathrm{D^b}(Y) \otimes_{\mathrm{D^bCoh}(\PP^1)} \mathrm{D^bCoh}(p)$ for $p \in \PP^1$. 
Even in this case, we expect to construct a mirror non-commutative Landau--Ginzburg model of $X_i$ 
by splitting the non-commutative Calabi--Yau fibration $W$ in an appropriate sense. 
It is not clear how to formulate a non-commutative Landau--Ginzburg model at this point, and defining a suitable formulation is part of the problem.

%%%%%%%%%%%%%%%%%%%%%%%%%%%%%%%%%%%%%%%%%%%%%%%%%%%%%%%%%%%%%%%%%%%%%%%%%%%%%%%%%%%%%%%%%%%%%
%%%%%%%%%%%%%%%%%%%%%%%%%%%%%%%%%%%%%%%%%%%%%%%%%%%%%%%%%%%%%%%%%%%%%%%%%%%%%%%%%%%%%%%%%%%%%

\section{Geometric quantization and theta functions} \label{GQ}
The aim of this section is to provide some evidence for the appearance of theta functions in mirror symmetry, for example, the DHT conjecture. 
A standard way to understand this is given by homological mirror symmetry, but we take a different path, namely geometric quantization. 
This section is somewhat disparate nature, compared with the previous sections, but we include it here because the subject is lurking and only partially explored in the context of mirror symmetry.

%%%%%%%%%%%%%%%%%%%%%%%%%%%%%%%%%%%%%%%%%%%%%%%%%%%%%%%%%%%%%%%%%%%%%%%%%%%%%%%%%%%%%%%%%%%%%
%%%%%%%%%%%%%%%%%%%%%%%%%%%%%%%%%%%%%%%%%%%%%%%%%%%%%%%%%%%%%%%%%%%%%%%%%%%%%%%%%%%%%%%%%%%%%

\subsection{Canonical basis via homological mirror symmetry} \label{basis HMS}
The appearance of theta functions in mirror symmetry is well-observed.  
Classically theta functions show up as a canonical basis of the vector space $H^0(A,L)$ for an ample line bundle $L$ of an abelian variety $A$. 
The definition of theta functions depends crucially on the group law of $A$, leaving us the impression that they are special to the abelian varieties. 
In the development of our understanding of mirror symmetry we come to expect that such a canonical basis exists for a much larger class of varieties. 

Heuristics from mirror symmetry is useful.  
According to Kontsevich's homological mirror symmetry \cite{Kon}, we have an equivalence of triangulated categories
$$
\mathrm{D^bCoh}(X)\cong \mathrm{D^bFuk}(Y).
$$ 
This implies an isomorphism between morphism spaces 
$$
Ext^*(E_1,E_2) \cong HF^*(L_1,L_2), 
$$
where $E_i$ and $L_i$ are mirror objects. 
Under a suitable condition, the Floer homology $HF^*(L_1,L_2)$ has a canonical basis (or set of generators) given by the intersections $L_1 \cap L_2$. 
Therefore homological mirror symmetry claims that $Ext^*(E_1,E_2)$ also has a canonical basis, which are often called {\it theta functions}. 
A prototypical example is given by an abelian variety $X$. 
If we take $E_1$ and $E_2$ to be the structure sheaf $\mathcal{O}_X$ and an ample line bundle $L$ respectively, 
then $Ext^*(E_1,E_2)=H^0(X,L)$ by the Kodaira vanishing theorem. 
In this case the canonical basis is nothing but the classical theta functions with characteristics.   
Since the mirror object of the sky scraper sheaf $\mathcal{O}_x$ is a Lagrangian torus $T \subset Y$, which gives rise to an SYZ fibration $\pi:Y\rightarrow B$, 
the above mirror correspondence and thus the canonical basis really depend upon a choice of an SYZ fibration. 

We elucidate the above theta functions for Calabi--Yau manifolds by the following examples. 
\begin{Ex}
Let us consider the Hasse pencil of cubic curves
$$
\{C_t\}=\{ x_0^3+x_1^3+x_2^3+3 t x_0x_1x_2=0 \} \subset \PP^2 \times \PP^1. 
$$
Then $\{x_i\}_{i=0}^3$ are the theta functions of a smooth Hasse cubic curve $C_t$.   
This is justified by observing that a smooth Hasse cubic curve is precisely the image of a $(3)$-polarized elliptic curve $C=\C /\Z +\tau \Z$ via the embedding 
$$
C \hookrightarrow \PP^2, \ \zeta \mapsto [\vartheta_{0,0}(3\zeta,3\tau),\vartheta_{\frac{1}{3},0}(3\zeta,3\tau),\vartheta_{\frac{2}{3},0}(3\zeta,3\tau)].
$$ 
In a similar manner, we can consider the Dwork pencil of quartic K3 surfaces
$$
\{Q_t\}=\{ x_0^4+x_1^4+x_2^4+x_3^4 + 4 t x_0x_1x_2x_3=0 \} \subset \PP^3 \times \PP^1. 
$$
Then $\{x_i\}_{i=0}^3$ are the theta functions of a smooth quartic K3 surface $Q_t$. 
\end{Ex}

In the following we however take a different path to observing theta functions, namely via geometric quantizations in mirror symmetry. 
We believe that quantization, degenerations and SYZ fibrations are the keys to the study of mirror symmetry. 

%%%%%%%%%%%%%%%%%%%%%%%%%%%%%%%%%%%%%%%%%%%%%%%%%%%%%%%%%%%%%%%%%%%%%%%%%%%%%%%%%%%%%%%%%%%%%
%%%%%%%%%%%%%%%%%%%%%%%%%%%%%%%%%%%%%%%%%%%%%%%%%%%%%%%%%%%%%%%%%%%%%%%%%%%%%%%%%%%%%%%%%%%%%

\subsection{Geometric quantization}
Geometric quantization is a recipe to construct a quantum theory out of a symplectic manifold $(X,\omega)$, which we think of as a classical phase space. 
Roughly we construct various objects associated to $(X, \omega)$ in such a way that certain analogues between the classical theory and the quantum theory remain manifest: 
\begin{itemize}
\item symplectic manifold $(X, \omega)$ $\rightsquigarrow$ Hilbert space $\mathcal{H}$.
\item algebra $C^\infty(X)$ of smooth functions $f$ $\rightsquigarrow$ algebra $\mathcal{A}_\hbar$ of operators $\widehat{f}$ on $\mathcal{H}$
\item Lagrangian submanifold $M \subset X$ $\rightsquigarrow$ state vector $\phi_M \in \mathcal{H}$
\item symplectomorphism $\mathrm{Symp}(X)$ $\rightsquigarrow$ automorphism $\mathrm{Aut}(\mathcal{A}_\hbar)$
\item $\dots$
\end{itemize}
satisfying various conditions. 
For instance the assignment 
$$
C^\infty(X) \rightarrow \mathcal{A}_\hbar: f \mapsto \widehat{f}
$$
should be a Lie algebra homomorphism. 
Here $C^\infty(X)$ is endowed with the Poisson bracket and $\mathcal{A}_\hbar$ with the commutator as Lie brackets. 

We begin with the classical example $X=T^*\R^n$, the symplectic vector space with the symplectic structure 
$$
\omega=\frac{1}{\hbar} \sum_{i=1}^n dq^i \wedge dp_i,
$$ 
Here $\{q^i\}$ are the coordinate of $\R^n$, $\{p_i\}$ are the canonical coordinate of the fiber direction, and $\hbar$ is the Plank constant. 
Let $\mathcal{H}=L^2(\R^n)$ be the space of $L^2$-functions on $\R^n$ with respect to the Lebesgue measure. 
From the requirement that the map $C^\infty(X) \rightarrow \mathcal{A}_\hbar$ is a Lie algebra homomorphism, we obtain the canonical commutation relation 
$$
[\hat{q}^i,\hat{p}_i]=\sqrt{-1} \hbar \delta^i_j. 
$$
The standard quantization is given by $\hat{q}^i=- \sqrt{-1}q^i, \hat{p}_i=\hbar \frac{\partial}{\partial q_i}$.

Keeping this example in mind, we shall study a symplectic manifold $X$ equipped with an integral symplectic class $[\omega] \in H^2(X,\Z)$. 
Then there exists a Hermitian line bundle $L \rightarrow X$ with a unitary connection $\nabla$ such that $c_1(L,\nabla)=\omega$. 
Such a pair $(L,\nabla)$ is called a pre-quantum bundle of the symplectic manifold $(X,\omega)$. 
We consider (the $L^2$-completion of) the space $\mathcal{H}=\Gamma^\infty(X,L)$ of smooth sections as a Hilbert space. 
We define $\widehat{f} \in \mathrm{End}(\mathcal{H})$ by the formula 
$$
\widehat{f}(s)=\nabla_{V_f}s-\sqrt{-1}fs,
$$
where $V_f$ is the Hamiltonian vector field associated to $f$. 
We can take $k\omega$ as a symplectic form and in this case $L^{\otimes k}$ is the pre-quantum bundle. 
Then $\frac{1}{k}$ plays the same role as the Plank constant $\hbar$. 

%%%%%%%%%%%%%%%%%%%%%%%%%%%%%%%%%%%%%%%%%%%%%%%%%%%%%%%%%%%%%%%%%%%%%%%%%%%%%%%%%%%%%%%%%%%%%
%%%%%%%%%%%%%%%%%%%%%%%%%%%%%%%%%%%%%%%%%%%%%%%%%%%%%%%%%%%%%%%%%%%%%%%%%%%%%%%%%%%%%%%%%%%%%

\subsection{Polarizations and polarized sections}
However it is known that the above Hilbert space $\mathcal{H}=\Gamma^\infty(X,L)$ is too big to capture the actual physical theory. 
In this sense, $\mathcal{H}$ is often called a pre-quantum Hilbert space. 
To get a reasonable theory, we need to choose a Poisson-commuting set of $n$ variables on the $2n$-dimensional phase space. 
Then consider the sections which depend only on these chosen variables.  
This is done by choosing a polarization, which is a coordinate-independent description of such a choice of $n$ Poisson-commuting sections.  

Let $X$ be a symplectic $2n$-fold equipped with an integral symplectic class: $[\omega] \in H^2(X,\Z)$ as above. 
We extend $\omega$ on $TX\otimes \C$ complex linearly. 
\begin{Def}
A polarization $P$ is an integrable Lagrangian subbundle in $TX \otimes \C$. 
In other words, it is a complex subbundle $P\subset TX \otimes \C$ of rank $n$ satisfying $[P,P] \subset P$ and $\omega|_P=0$. 
\end{Def}

For a polarization $P$, we would like to define the space of polarized sections by 
$$
\Gamma_P(X,L)=\{ s \in \Gamma(X,L) \ | \ \nabla_v s =0 \ \forall v \in P\}. 
$$
This means the sections are constant in the $P$-direction. 
%We want the sections to be functions of only $n$ variables on the $2n$-dimensional classical phase space. 

It is not clear at this point how we can universally define polarized sections, and we may need some modification of this definition depending on a polarization. 
It is worth mentioning that the space of polarized sections is the space of wave functions in physics. 

There are many examples of polarizations but we will focus on the following two polarizations relevant to mirror symmetry: 
Lagrangian polarization and K\"ahler polarization. 

Given a Lagrangian fibration $\pi:X \rightarrow B$, then the complexified relative tangent bundle 
$$
P=T_{X/B}\otimes \C =\Ker(d\pi: TX \rightarrow TB)\otimes \C
$$
gives a polarization, which we call a Lagrangian polarization.  
This polarization is real in the sense that $\overline{P}=P$. 

Let us assume further that $(X,\omega)$ be a K\"ahler manifold. 
Then the anti-holomorphic tangent bundle 
$$
P=T^{0,1}X \subset TX\otimes \C
$$
is a polarization. Such a polarization is called a K\"ahler polarization. 
Note that $L$ is holomorphic since the curvature $\omega$ is of type $(1,1)$. 
A K\"ahler polarization is characterized by the conditions $P \cap \overline{P}=0$ and $\omega|_{P\times \overline{P}}>0$. 
Giving a K\"ahler polarization is equivalent to fixing a compatible complex structure on the underlying symplectic manifold $(X,\omega)$. 
In this case, the space of polarized sections $\Gamma_P(X,L)$ is nothing but the space of holomorphic sections $H^0(X,L)$.

%%%%%%%%%%%%%%%%%%%%%%%%%%%%%%%%%%%%%%%%%%%%%%%%%%%%%%%%%%%%%%%%%%%%%%%%%%%%%%%%%%%%%%%%%%%%%
%%%%%%%%%%%%%%%%%%%%%%%%%%%%%%%%%%%%%%%%%%%%%%%%%%%%%%%%%%%%%%%%%%%%%%%%%%%%%%%%%%%%%%%%%%%%%

\subsection{BS Lagrangian submanifolds and quantization problem}
Let $X$ be a symplectic $2n$-fold as above and $M \subset X$ a Lagrangian submanifold. 
Since $\omega$ is the curvature of $M$, the restriction $L|_M$ is a flat line bundle on $M$.  

\begin{Def}
A Lagrangian submanifold $M \subset X$ is called Bohr--Sommerfeld if the restriction $(L,\nabla)|_M$ is trivial. 
More generally, $M$ is called Bohr--Sommerfeld of level $k \in \N$ if the restriction $(L^{\otimes k},\nabla)|_{M}$ is trivial. 
\end{Def}

Geometrically, a Lagrangian submanifold $M\subset X$ is Bohr--Sommerfeld of level $k$ if and only if $k\int_D \omega \in \Z$ for every disk $D \subset X$ such that $\partial D \subset M$.  

Let $\pi:X\rightarrow B$ be a Lagrangian torus fibration. 
We assume that the image $\pi(X)$ of the moment map is an integral polytope in $\R^n$ (this can always be realized by a suitable shift).  
Then a fiber $\pi^{-1}(b)$ is a Bohr--Sommerfeld Lagrangian submanifold of level $k$ if and only if $b \in \frac{1}{k}\Z^n\cap \Delta \subset \R^n$.  
In particular, Bohr--Sommerfeld Lagrangian fibers appear discretely. 

As we discussed in Section \ref{sLag}, the special Lagrangian submanifolds form a nice class of Lagrangian submanifolds.  
An advantage of the Bohr--Sommerfeld Lagrangians is that they are defined purely in symplectic geometry, no need for a Calabi--Yau structure. 

\begin{Ex} \label{toric BS}
Let $\PP_\Delta$ be the toric $n$-fold associated to a lattice polytope $\Delta \subset \R^n$.  
It naturally comes with a polarization $L$ and we consider the toric K\"ahler form $\omega=c_1(L)$. 
Then a torus fiber $\pi^{-1}(b)$ of the moment map $\pi_\Delta:\PP_\Delta \rightarrow \Delta$ is a Bohr--Sommerfeld Lagrangian torus of level $k$ if and only if $b \in \frac{1}{k}\Z^n\cap \mathrm{int}(\Delta)$. 
It makes sense to speak about Bohr--Sommerfeld isotropic submanifolds. 
If we think of isotropic fibers as degenerate Lagrangian submanifolds, then we observe that the possibly degenerate Bohr--Sommerfeld Lagrangian fibers are indexed by $\frac{1}{k}\Z^n\cap \Delta$. 
\end{Ex}

There is in general no non-trivial section of $L$ which is constant along the fibers because $L|_{\pi^{-1}(b)}$ may have non-trivial holonomy. 
%This means that the space of Lagrangian polarized smooth sections is typically trivial. 
Hence instead of the smooth sections, it seems reasonable to consider the distributional sections supported on the Bohr--Sommerfeld fibers. 
So we shall modify the space of polarized sections as follows: 
\begin{align}
\Gamma_{X/B}(X,L)&=\{ s \in \Gamma^{\mathrm{dist}}(X,L) \ | \ \mathrm{supp}(s) \subset \mathrm{BS \ fibers}, \  \nabla_v s =0 \ \forall v \in P\}. \notag \\
&=\oplus \C \delta_{BSL}. \notag 
\end{align}
where $\delta_{BSL}$ denotes the delta function on a possibly degenerate Bohr--Sommerfeld Lagrangian fiber. 
The last equality follows from the fact that covariantly constant sections on a Bohr--Sommerfeld Lagrangian fiber is unique up to multiplication by constants. 
Therefore the Bohr--Sommerfeld Lagrangian fibers form a canonical basis of the space of Lagrangian polarized sections as a vector space.

Let $L \rightarrow X$ be a pre-quantum bundle for a compact K\"ahler manifold $X$, or equivalently a polarized K\"ahler manifold $(X,L)$. 
Given a Lagrangian torus fibration $\pi:X\rightarrow B$,  
then we have two space of wave functions $\Gamma_{T_{M/B}}(X,L)$ and $H^0(X,L)$.

\begin{Conj}[Quantization Problem] \label{QP}
These two spaces of wave functions are isomorphic (in a canonical way): $$\Gamma_{X/B}(X,L) \cong H^0(X,L).$$ 
\end{Conj}
The basic philosophy of the conjecture is that there may be several ways to quantize a given physical system but the end result should not depend on a way of quantization. 

There are several important examples where this conjecture is shown to be true. 
For instance, Guillemin and Sternberg proved the conjecture for the flag manifolds viewed as the Gelfand--Cetlin integrable system \cite{GuiSte}. 
They did not give any direct relationship between the quantizations, 
but recently Hamilton and Konnno \cite{HamKon} described a deformation of the complex structure on the flag manifold leading to a convergence of polarizations joining the two quantizations.
Another important result is due to Anderson, who showed $$\Gamma_{X/B}(X,L^{\otimes k}) \cong H^0(X,L^{\otimes k})$$ for sufficiently large $k$ 
provided that the Lagrangian fibration $\pi:X\rightarrow B$ has no degenerate fiber \cite{And}. 
It is of interest to observe that the conjecture holds for a sufficiently small $k$ or {\it a sufficiently small scale}.  
In general it is very difficult to compute $\Gamma_{X/B}(X,L)$ as there exist degenerate fibers. 

One important implication of Conjecture \ref{QP} is that $H^0(X,L)$ has a canonical basis corresponding to $\{\delta_{BSL}\}\subset \Gamma_{X/B}(X,L)$. 
Recall that choosing a (special) Lagrangian torus fibration of a Calabi--Yau manifold is equivalent to choosing a mirror manifold $Y$. 
Then $H^0(X,L)$ has two canonical bases: one induced by homological mirror symmetry and one induced by geometric quantization.  
Then a natural question is whether they coincide or not.

\begin{Conj}
Let $(X,L)$ be a polarized Calabi--Yau manifold.  
Given a (special) Lagrangian torus fibration $\phi:X\rightarrow B$, then $H^0(X,L^{\otimes k})$ has a canonical basis given in Theorem \ref{QP} for a sufficiently large $k$. 
Moreover, it corresponds to the canonical basis of the Floer cohomology under the mirror correspondence associated to $\phi$.  
\end{Conj}

\begin{Ex}[Example \ref{toric BS} continued] \label{toric monomial} 
For the toric moment map $\pi_\Delta:\PP_\Delta \rightarrow \Delta$, the space of Lagrangian polarized sections $\Gamma_{\PP_\Delta/\Delta}(\PP_\Delta,L)$ has a basis indexed by $\frac{1}{k}\Z^n\cap \Delta$.  
On the other hand, it is well-know in toric geometry that $ H^0(\PP_\Delta,L)$ has a canonical basis, called the monomial basis indexed by $\frac{1}{k}\Z^n\cap \Delta$. 
Therefore we observe that there is a natural isomorphism $$\Gamma_{\PP_\Delta/\Delta}(\PP_\Delta,L) \cong H^0(\PP_\Delta,L)$$ given by the identification of the two canonical bases. 
\end{Ex}

In fact, in the toric case, more precise correspondence is given by Baier, Florentino, Mourao, and Nunes in \cite{BFMN}. 
They carried out, by changing symplectic potentials, a deformation of toric K\"ahler polarization which joins the two poalrizations. 
This deformation explicitly shows that the monomial basis converge to delta-function sections supported on the Bohr--Sommerfeld Lagrangian fibers. 

A naive way to construct such a basis for a Calabi--Yau manifold $X$ is again to consider a toric degeneration $X \rightsquigarrow \cup_{\Delta}\PP_\Delta$. 
As Example \ref{toric monomial} illustrates, the conjecture holds for the union of toric varieties $\cup_{\Delta}\PP_\Delta$. 
Thus it is reasonable to expect that Bohr--Sommerfeld Lagrangian fibers of $\cup_\Delta \pi_\Delta$ deform to those of $\pi$, 
where $\pi:X\rightarrow B$ is a Lagrangian torus fibration approximated by the toric degeneration (Section \ref{approximate SYZ}). 

\subsection{BS Lagrangian submanifolds and mirror symmetry}

Finally let us try to understand Bohr--Sommerfeld Lagrangian submanifolds from the viewpoint of mirror symmetry. 
First recall the semi-flat mirror symmetry associated to an integral affine manifold $B$. 
$$ 
\xymatrix{
TB/\Lambda \ar[rd]_\phi &  & T^*B/\Lambda^* \ar[ld]^{\phi^\vee}\\
 & B & 
}
$$
We take a local coordinate $(z_i)=\exp(- 2\pi(x_i+\sqrt{-1}y_i))$ of $TB/\Lambda$. 
Recall that, for $b \in B$,  a point on a Lagrangian fiber $(\phi^\vee)^{-1}(b)$ corresponds to a flat $\U(1)$-connection on the dual fiber $L_b=\phi^{-1}(b)$. 
Therefore the real Fourier--Mukai transform maps the graph $\Gamma(s)$ of a section $s=\sum{i}s_i(x) dx_i$ of $T^*B/\Lambda^* \rightarrow B$ to 
the trivial bundle $L$ on $TB/\Lambda$ with the connection 
$$
\nabla=d+2 \pi \sqrt{-1}\sum_i s_i(x) dy_i. 
$$
The $(2,0)$-part of the curvature of $\nabla$ is given by 
$$
F^{2,0}_\nabla=\pi \sum_{i,j}(\frac{\partial s_j}{\partial x_i}-\frac{\partial s_i}{\partial x_j})dz_i \wedge dz_j.
$$
Then it follows that the graph $\Gamma(s)$ is Lagrangian (i.e. $s$ is closed) if and only $\nabla$ defines a holomorphic line bundle on $TB/\Lambda$. 
Moreover, $\nabla|_{L_b}$ is trivial if and only if $b \in \Gamma(s) \cap \Gamma(0)$. 
In particular, if $L$ defines an ample line bundle on $TB/\Lambda$ and $TB/\Lambda$ is considered as a symplectic manifold with respect to $c_1(L)$, 
the fibers $\{L_b\}_{b  \in \Gamma(s) \cap \Gamma(0)}$ are precisely the Bohr--Sommerfeld Lagrangian fibers. 

Now let us look at the Tyurin degeneration of an elliptic curve $X\rightsquigarrow X_1 \cup_Z X_2$ (Section \ref{DHT proof}).  
In this case, we have a natural pair of Lagrangian sections of a Lagrangian fibration. 
One is a chosen zero section $s_1$, the other is the image $s_2$ of the zero section under the Picard--Lefschetz monodromy. 
Then the monodromy matrix with respect to the zero section and the vanishing cycle is given by $\begin{bmatrix} 1 & 2 \\  0 & 1 \\ \end{bmatrix}$. 
Therefore the corresponding fibration $Y \rightarrow \PP^1$ on the mirror side is induced 
by the $(2)$-polarization (recall the duality between $H^{1,0}(X)\rightarrow H^{0,1}(X)$ and $H^{0,0}(Y)\rightarrow H^{1,1}(Y)$). 
Hence we observe a compatible fact that $\Gamma(s_1) \cap \Gamma(s_2)$ consists of 2 points and $(2)$-polarization has 2 theta functions as a basis of the space of its global sections. 

Despite the fundamental nature of the Bohr--Sommerfeld Lagrangian submanifolds, many basic questions are still unanswered. 
Mirror symmetry has already sparked many exiting development of mathematics and the author wishes to continue to see many more of such development in geometric quantization. 

%%%%%%%%%%%%%%%%%%%%%%%%%%%%%%%%%%%%%%%%%%%%%%%%%%%%%%%%%%%%%%%%%%%%%%%%%%%%%%%%%%%%%%%%%%%%%
%%%%%%%%%%%%%%%%%%%%%%%%%%%%%%%%%%%%%%%%%%%%%%%%%%%%%%%%%%%%%%%%%%%%%%%%%%%%%%%%%%%%%%%%%%%%%

%%%%%%%%%%%%%%%%%%%%%%%%%%%%%%%%%%%%%%%%%%%%%%%%%%%%%%%%%%%%%%%%%%%%%%%%%%%%%%%%%%%%%%%%%%%%%
%%%%%%%%%%%%%%%%%%%%%%%%%%%%%%%%%%%%%%%%%%%%%%%%%%%%%%%%%%%%%%%%%%%%%%%%%%%%%%%%%%%%%%%%%%%%%

\par\noindent{\scshape \small
Department of Mathematics, Kyoto University\\
Kitashirakawa-Oiwake, Sakyo, Kyoto, 606-8502, Japan}
\par\noindent{\ttfamily akanazawa@math.kyoto-u.ac.jp}
\end{document}